\documentclass{amsart}
\usepackage[margin=1in]{geometry}
\usepackage{shor-research}
\usepackage{diagbox}
\usepackage{caption}
\newcommand{\Ceil}[1]{\left\lceil #1 \right\rceil} 
\newcommand{\fpf}[2]{\left\{\frac{#1}{#2}\right\}} 
\newcommand{\flrf}[2]{\left\lfloor\frac{#1}{#2}\right\rfloor} 
\newcommand{\clf}[2]{\left\lceil\frac{#1}{#2}\right\rceil} 
\newcommand{\divs}{\mathrm{d}} 
\newcommand{\Jacr}{\mathrm{r}} 
\newcommand{\DirichletD}{\mathrm{D}}
\newcommand{\GaussC}{\mathrm{C}}
\newcommand{\ee}{e}
\newcommand{\ii}{i}
\newcommand{\dx}{\mathrm{d}x}

\newcommand{\R}{\mathbb{R}}

\newcommand{\Fseq}{\mathcal{F}} 
\newcommand{\Cseq}{\mathcal{C}} 
\newcommand{\Rseq}{\mathcal{R}} 
\newcommand{\Num}{\mathcal{N}} 
\newcommand{\reallydivs}{\dagger} 
\newcommand{\genlegendre}[4]{%
  \genfrac{(}{)}{}{#1}{#3}{#4}%
  \if\relax\detokenize{#2}\relax\else_{\!#2}\fi
}

\title{On residues of rounded shifted fractions with a common numerator}
\author{Nicholas Dent \and Caleb M.\ Shor}
\date{\today}

\DeclareMathOperator{\nint}{nint}
\begin{document}
\maketitle
\begin{abstract}
    For any positive integer $n$ along with parameters $\alpha$ and $\nu$, we define and investigate $\alpha$-shifted, $\nu$-offset, floor sequences of length $n$. We find exact and asymptotic formulas for the number of integers in such a sequence that are in a particular congruence class. As we will see, these quantities are related to certain problems of counting lattice points contained in regions of the plane bounded by conic sections. We give specific examples for the number of lattice points contained in elliptical regions and make connections to a few well-known rings of integers, including the Gaussian integers and Eisenstein integers.
\end{abstract}
\section{Introduction}
For a fixed positive integer $n$, consider the integer sequence 
\begin{equation}\label{seq:intro}
    \Floor{\frac{n}{1}}, \Floor{\frac{n}{2}}, \Floor{\frac{n}{3}}, \dots, \Floor{\frac{n}{n}},
\end{equation}
where $\Floor{x}$ denotes the floor function for real $x$. Among the $n$ terms in this sequence, how many are odd? At first glance, it maybe seems reasonable to expect the proportion of odd numbers in the sequence to be roughly half. For $n=10$, the sequence is $10, 5, 3, 2, 2, 1, 1, 1, 1, 1$, of which $7/10=70\%$ are odd.  In general, terms in the second half of the sequence are all 1 and thus odd, implying the proportion of odd terms is always at least $50\%$. Looking at more data, this proportion appears to hover around $69\%$ as $n$ grows.

This leads to a few natural questions. Why $69\%$? What happens if we replace the floor function with the ceiling function? Or if we round to the nearest integer?

\begin{figure}
    \centering
    \includegraphics[width=.6\textwidth]{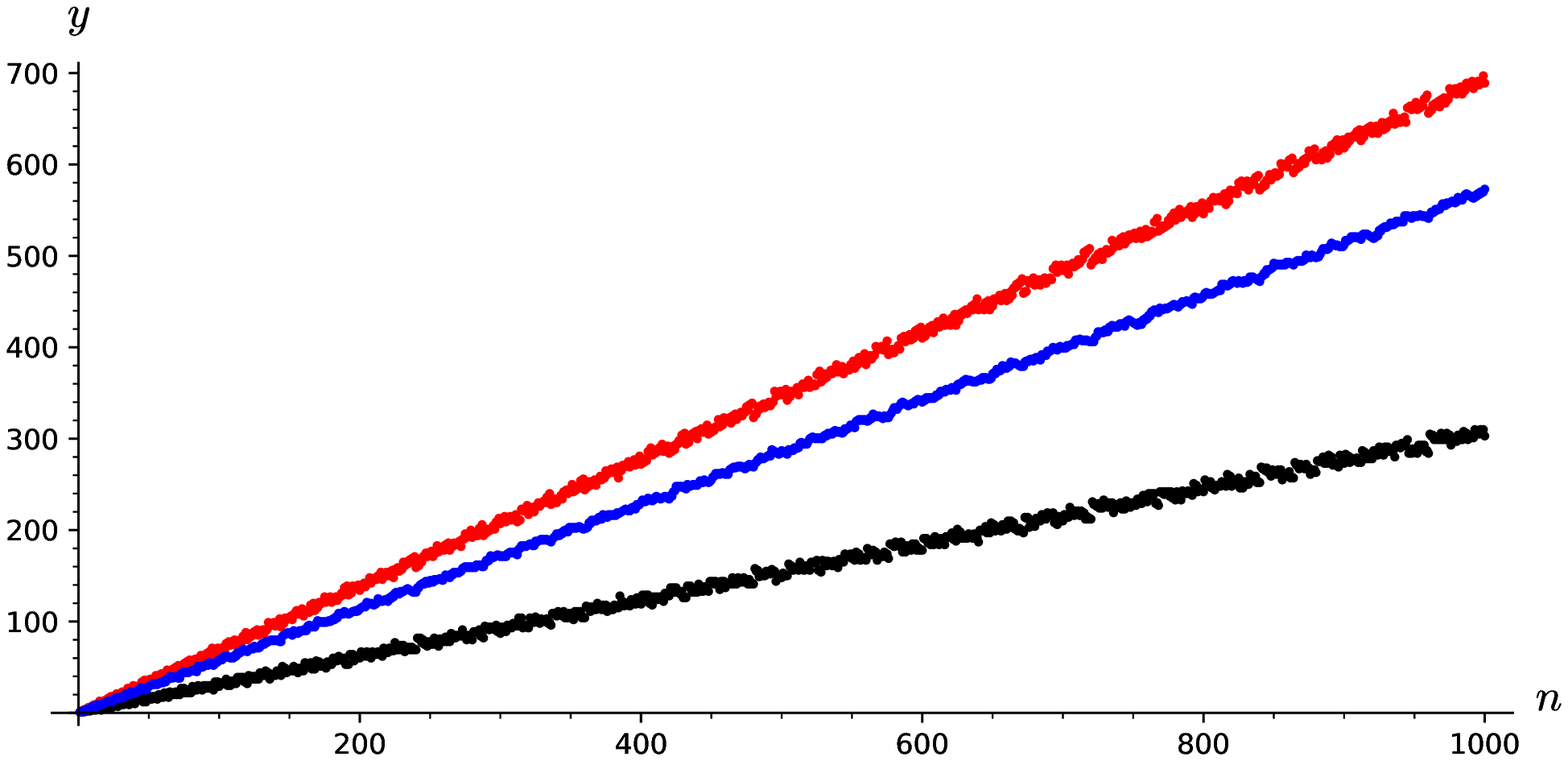}
    \caption{Plots of $y=\Fseq_n$ (red, top-most), $y=\Cseq_n$ (blue, middle), and $y=\Rseq_n$ (black, bottom-most), for $1\le n\le 1000$}
    \label{fig:first-three-sequences}
\end{figure}

In order to answer these questions, we consider three integer sequences defined, for positive integers $n$, in terms of the floor, ceiling, and nearest integer rounding functions:
\begin{align*}
    \Fseq_n &= \#\left\{k\in\Z : 1\le k\le n,\, \Floor{n/k}\text{ is odd}\right\},\\
    \Cseq_n &= \#\left\{k\in\Z : 1\le k\le n,\, \Ceil{n/k}\text{ is odd}\right\},\\
    \Rseq_n &= \#\left\{k\in\Z : 1\le k\le n,\, \nint(n/k)\text{ is odd}\right\},
\end{align*}
where $\Ceil{x}$ denotes the ceiling function and $\nint(x)$ the nearest integer rounding function\footnote{The definition of the nearest integer function can be ambiguous for half-integers. A common convention is to round to the nearest even integer. Since we are interested in parity, we instead choose to always round half-integers up. Hence, $\nint(2.5)=3$, $\nint(3.5)=4$, and so on. As we will see, the asymptotic formula for $\Rseq_n$ will be the same regardless of how we choose to round half-integers.}.
\begin{figure}
    \centering
    \begin{tabular}{c|rrrrrrrrrrrrrrrrrrrr}
    $n$ & 1 & 2 & 3 & 4 & 5 & 6 & 7 & 8 & 9 & 10 & 11 & 12 & 13 & 14 & 15 & 16 & 17 & 18 & 19 & 20 \\ \hline 
    $\Fseq_n$ & 1 & 1 & 3 & 2 & 4 & 4 & 6 & 4 & 7 & 7 & 9 & 7 & 9 & 9 & 13 & 10 & 12 & 12 & 14 & 12 \\ 
    $\Cseq_n$ & 1 & 1 & 2 & 1 & 3 & 2 & 3 & 2 & 5 & 3 & 4 & 3 & 6 & 5 & 6 & 3 & 7 & 6 & 7 & 6 \\ 
    $\Rseq_n$ & 1 & 1 & 2 & 2 & 4 & 3 & 4 & 4 & 6 & 7 & 6 & 5 & 9 & 8 & 9 & 9 & 10 & 10 & 11 & 12
\end{tabular}
    \caption{Terms in the sequences $\Fseq_n$, $\Cseq_n$, $\Rseq_n$ for $1\le n\le 20$}
    \label{fig:first-20-terms}
\end{figure}
The first 20 terms of each of these three sequences can be found in Figure~\ref{fig:first-20-terms}. The sequences $\Fseq_n$ and $\Cseq_n$ appear in The On-Line Encyclopedia of Integer Sequences (\cite{OEIS}) as, respectively, sequences \href{https://oeis.org/A059851}{A059851} and \href{https://oeis.org/A330926}{A330926}. 
Plots of $\Fseq_n$, $\Cseq_n$, and $\Rseq_n$, for $1\le n\le 1000$, appear in Figure~\ref{fig:first-three-sequences}.

The plots in Figure~\ref{fig:first-three-sequences} suggest that each sequence is roughly linear. As we will show, this is true asymptotically. We will find asymptotic formulas for each of these sequences, and from those we will conclude that
\[
    \lim\limits_{n\to\infty}\dfrac{\Fseq_n}{n}=\log2,
    \quad
    \lim\limits_{n\to\infty}\dfrac{\Cseq_n}{n}=1-\log2,
    \quad\text{ and } 
    \lim\limits_{n\to\infty}\dfrac{\Rseq_n}{n}=\frac{\pi}{2}-1.
\]
In particular, the proportion of odd terms in $\Fseq_n$ is asymptotically $\log2\approx 0.693147$. This explains the $69\%$ mentioned in the first paragraph.

Our results follow from two classical results in analytic number theory: the \emph{Dirichlet divisor problem} (Theorem~\ref{thm:dirichlet}), in which one wants to count the number of lattice points in the first quadrant of the plane beneath a hyperbola; and the \emph{Gauss circle problem} (Theorem~\ref{thm:gauss}), in which one wants to count the number of lattice points contained within a circle in the plane. 

Each of these classical results leads to a geometric interpretation for one of our sequences. In Proposition~\ref{prop:F_n-Dirichlet}, we will see that $\Fseq_n$ is the difference of the numbers of lattice points in two hyperbolic regions in the plane. In Proposition~\ref{prop:gauss-R_n}, we will see that the number of lattice points contained in a circle of radius $\sqrt{2n}$ is equal to $4\Rseq_n+4n+1$.

We then consider more general sequences. Since $\nint(x)=\Floor{x+1/2}$, we can think of the nearest integer function as the usual floor function shifted by $1/2$. If we replace this $1/2$ with an arbitrary real number $\alpha$ and also incorporate a real number $\nu$ to offset the numerator $n$, we obtain an \emph{$\alpha$-shifted, $\nu$-offset, floor sequence of length $n$}
\begin{equation}\label{seq:alpha-shift-intro}
    \Floor{\frac{n-\nu}{1}+\alpha}, 
    \Floor{\frac{n-\nu}{2}+\alpha}, 
    \dots, 
    \Floor{\frac{n-\nu}{n}+\alpha}.
\end{equation}
We can determine how many of the integers in Sequence~\eqref{seq:alpha-shift-intro} are odd. A natural problem, which we will address, is to find $\alpha$ and $\nu$ for which (asymptotically) half of these integers are odd and half are even. We will see that there is a unique value for $\alpha\in[0,1]$, independent of $\nu$, for which this occurs, and numerically approximate it.

Another problem is to count the number of integers in Sequence~\eqref{seq:alpha-shift-intro} that belong to a given congruence class of integers. We will find exact and asymptotic formulas for these counts. In particular, by Corollary~\ref{cor:N_m-slope}, for any $m\in\N$ and $\alpha\in[0,1)$, the proportion of integers in an $\alpha$-shifted, $\nu$-offset, floor sequence of length $n$ that are congruent to 1 modulo $m$ is asymptotically equal to 
\[
    \frac{-\alpha}{1-\alpha}+
    \int\limits_0^1 \frac{(1-x)x^{-\alpha}}{1-x^m}\dx
\]
and, for $2\le r\le m$, the proportion of such integers that are congruent to $r$ modulo $m$ is asymptotically equal to 
\[
    \int\limits_0^1 \frac{(1-x)x^{r-1-\alpha}}{1-x^m}\dx.
\]

Finally, we will look more closely at connections between the number of integers in a certain congruence class in an $\alpha$-shifted, $\nu$-offset, floor sequence of length $n$ and the number of lattice points in a certain region of the plane. We will demonstrate the connections by obtaining formulas in terms of these counts for the number of lattice points in the following elliptical regions: $x^2+y^2\le n$, $x^2+xy+y^2\le n$, and $x^2+2y^2\le n$, each for any $n\in\N$. Our results here follow from the theory of binary quadratic forms.

\subsection{Notation}
$\Z$ and $\N$ denote, respectively, the sets of integers and positive integers. All logarithms in this paper use base $\ee$. The cardinality of a finite set $A$ is denoted $\#A$. We use big O notation as follows. For functions $f(x)$ and $g(x)$, if there exist constants $M$ and $a$ such that $|f(x)|\le M g(x)$ for all $x\ge a$, then we write $f(x)=\bigO(g(x))$.

\subsection{Organization}
This paper is organized as follows. In Section~\ref{sec:F_n and C_n}, we find exact and asymptotic formulas for $\Fseq_n$ and $\Cseq_n$. In Section~\ref{sec:Rseq_n}, we do the same for $\Rseq_n$. In Section~\ref{sec:generalized}, for integers $r,m$ with $1\le r\le m$, and for $\alpha\in[0,1)$, we find exact and asymptotic formulas for $\Num_{n,\alpha,\nu,r,m}$, the number of integers in the $\alpha$-shifted, $\nu$-offset, floor sequence of length $n$ that are congruent to $r$ modulo $m$. Finally, in Section~\ref{sec:applications}, we have two tasks. First, we will compute the (unique) shift $\alpha=\alpha_0$ for which (asymptotically) the $\alpha$-shifted, $\nu$-offset, floor sequence of length $n$ contains as many odd terms as even terms. Second, we will determine the number of lattice points in certain elliptical regions of the plane in terms of the sequences $\Num_{n,\alpha,\nu,r,m}$.

\section{The floor and ceiling sequences}\label{sec:F_n and C_n}

\subsection{The floor sequence}
As was mentioned in the introduction, for $n\in\N$, $\Fseq_n$ is equal to the number of odd integers in the floor sequence 
\begin{equation}\label{seq:floor}
    \Floor{\frac{n}{1}}, \Floor{\frac{n}{2}}, \dots, \Floor{\frac{n}{n}},
\end{equation}
where the floor function is defined, for $x\in\R$, by $\Floor{x}=\max\{z\in\Z : z\le x\}$.  In this subsection, we will find an exact formula for $\Fseq_n$. Then, with that formula and the solution to Dirichlet's divisor problem, we will find an asymptotic formula for $\Fseq_n$.

We start by counting the number of integers in Sequence~\eqref{seq:floor} that are greater than or equal to a given integer.
\begin{lem}\label{lem:consec_flr}
    Fix $n\in\N$. For $k\in \N$,  Sequence~\eqref{seq:floor} 
    contains $\Floor{n/k}$ terms that are greater than or equal to $k$.
\end{lem} 

\begin{proof}
Via the division algorithm, there are $q,r\in \Z$ for which $n=kq+r$ with $0\leq r< k$. Then $\Floor{n/k}=q$. We will show that the first $q$ terms in Sequence~\eqref{seq:floor} are greater than or equal to $k$. 

Note that 
\[
    \flrf{n}{q} 
    = \flrf{kq+r}{q} 
    = \Floor{k+\frac{r}{q}}\ge k.
\] 
Now consider an arbitrary term $\Floor{n/d}$ in Sequence~\eqref{seq:floor}. Then $1\leq d\leq n$. If $d < q$, then $n/d>n/q$ and hence $\Floor{n/d}\ge\Floor{n/q}\ge k$. On the other hand, if $d > q$, then $d\ge q+1$. This means $kd \ge k(q+1) > n$ and therefore $k>n/d\geq \Floor{n/d}$. We have already shown that $\Floor{n/d}\geq k$ if $d=q$. We conclude that $\Floor{n/d}\geq k$ if and only if $d \in \left\{1,2,\dots ,q\right\}$.
\end{proof} 

We can now find an exact formula for $\Fseq_n$.
\begin{prop}\label{prop:floor_seq}
    For $n\in\N$, 
    \[
        \Fseq_n
        =\sum\limits_{d=1}^n\Floor{\frac{n}{d}}(-1)^{d+1}.
    \]
\end{prop}
\begin{proof}
We wish to count the number of odd terms in Sequence~\eqref{seq:floor}. Here we can use Lemma~\ref{lem:consec_flr} to see that for a given $k$ there are $\Floor{n/k}$ terms whose value is at least $k$, and thus $\Floor{n/k}-\Floor{n/(k+1)}$ terms whose value is exactly $k$. To compute $\Fseq_n$, we sum for all odd $k$. We find
    \[
        \Fseq_n 
        =\sum\limits_{\substack{1\le k\leq n \\ k \text{ odd}}} \left(\flrf{n}{k}-\flrf{n}{k+1}\right)
        =\sum\limits_{d=1}^n\Floor{\frac{n}{d}}(-1)^{d+1},
    \]
as desired.
\end{proof}

Next, we wish to find an asymptotic formula for $\Fseq_n$. To start, we will manipulate the formula for $\Fseq_n$ obtained in Proposition~\ref{prop:floor_seq} by considering odd and even $d$ separately.
\begin{align*}
    \Fseq_n
    &=\sum\limits_{d\text{ odd}} \Floor{\frac{n}{d}}-\sum\limits_{d\text{ even}} \Floor{\frac{n}{d}}\\
    &=\sum\limits_{d\text{ odd}} \Floor{\frac{n}{d}}-\sum\limits_{d\text{ even}} \Floor{\frac{n}{d}}+\sum\limits_{d\text{ even}} \Floor{\frac{n}{d}}-\sum\limits_{d\text{ even}} \Floor{\frac{n}{d}}\\
    &=\sum\limits_{d=1}^n\Floor{\frac{n}{d}} - 2\sum\limits_{b=1}^{\Floor{n/2}}\Floor{\frac{n}{2b}}\\
    &=\sum\limits_{d=1}^n\Floor{\frac{n}{d}} - 2\sum\limits_{b=1}^{\Floor{n/2}}\Floor{\frac{n/2}{b}}.
\end{align*}

We will now use $\tau(n)$, the multiplicative function which counts the number of positive divisors of $n$, and $\DirichletD(x)=\sum\limits_{m\le x}\tau(m)$, the \emph{divisor summatory function}, where the summation is taken over positive integers $m$. With this function, we have the following:
\[
    \DirichletD(n)
    = \sum\limits_{m=1}^n\tau(m)
    = \sum\limits_{m=1}^n\sum\limits_{d\mid m}1 
    = \sum\limits_{d=1}^n\Floor{\frac{n}{d}}.
\]
Geometrically, $\DirichletD(n)$ is the number of lattice points in the interior and boundary of the region in the $xy$-plane bounded by the graphs of $x=1$, $y=1$, and the hyperbola $xy=n$.

Observe that $\Fseq_n=\DirichletD(n)-2\DirichletD(n/2)$. For $n\in\N$, we define two regions. Let $H_{1,n}$ be the hyperbolic region in the $xy$-plane with $x,y\ge1$ bounded by the graphs of the hyperbola $xy=n$ and the hyperbola $xy=n/2$. Let $H_{2,n}$ be the region in the $xy$-plane bounded by the graphs of $x=1$, $y=1$, and the hyperbola $xy=n/2$. Then $\DirichletD(n)$ is equal to the number of lattice points in $H_{1,n}\cup H_{2,n}$, and $\DirichletD(n/2)$ is equal to the number of lattice points in $H_{2,n}$. This gives us a geometric interpretation for $\Fseq_n$.

\begin{prop}\label{prop:F_n-Dirichlet}
    For $n\in\N$, $\Fseq_n$ is equal to the number of lattice points in $H_{1,n}$ minus the number of lattice points in $H_{2,n}$. For $H_{1,n}$, we include all points on the boundary except for those on the boundary curve $xy=n/2$. For $H_{2,n}$, we include all points on the boundary.
\end{prop}
\begin{figure}
    \centering
    \includegraphics[scale=.7]{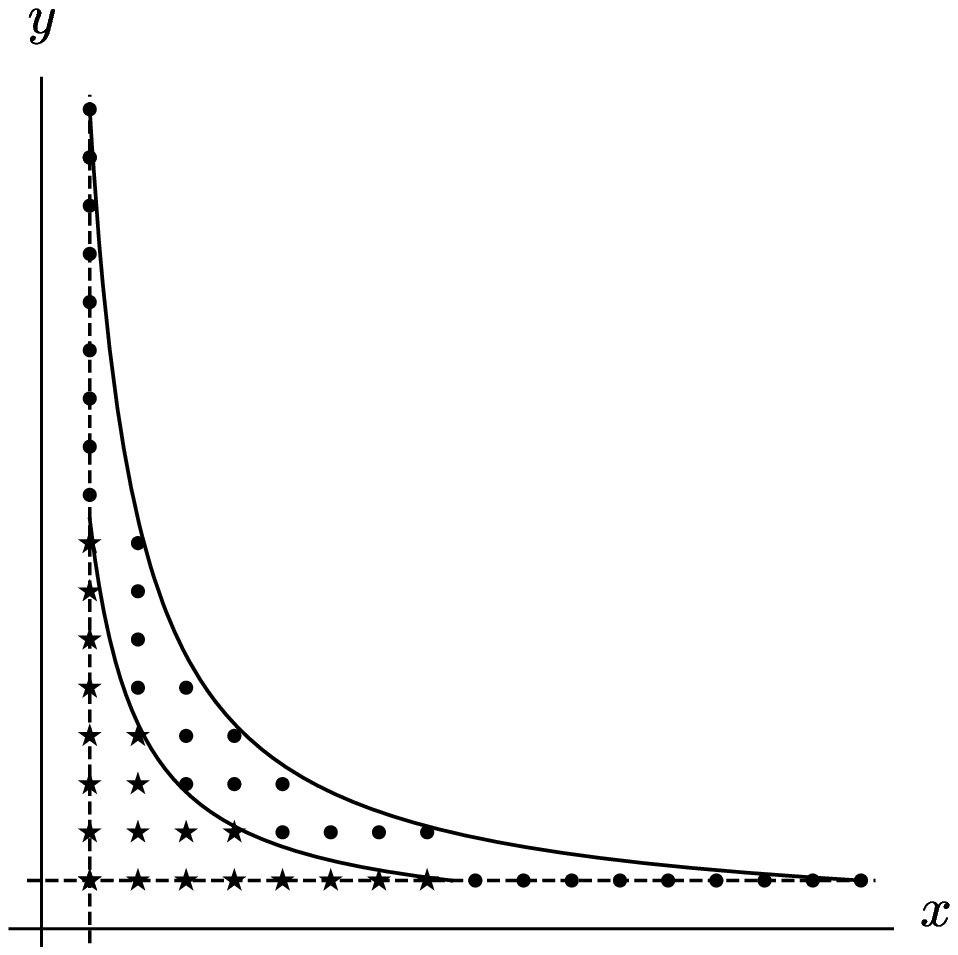}
    \caption{For $n=17$, the graphs of $xy=17$ and $xy=17/2$ along with lattice points in the hyperbolic regions $H_{1,17}$ (circles) and in $H_{2,17}$ (stars).}
    \label{fig:hyperbola-F_n}
\end{figure}
\begin{example}
    To illustrate Proposition~\ref{prop:F_n-Dirichlet}, we have drawn the regions $H_{1,n}$ and $H_{2,n}$ for $n=17$ in Figure~\ref{fig:hyperbola-F_n}. Since $H_{1,17}$ contains 32 points and $H_{2,17}$ contains 20 points, we find $\Fseq_{17}=32-20=12$.
\end{example}

If we have an asymptotic formula for $\DirichletD(x)$, known as the \emph{Dirichlet divisor problem}, then we get an asymptotic formula for $\Fseq_n$. The following theorem of Dirichlet is exactly what we need. (For details, see, e.g., \cite[Theorem 3.3]{Apostol1976}. We will use a slightly modified version of this approach in the proof of Proposition~\ref{prop:sum-of-floors-to-sum-of-fractions} later in this paper.)
\begin{thm}\label{thm:dirichlet}
    For all $x\ge1$, 
    \[
        \DirichletD(x)
        =x\log{x}+(2\gamma-1)x+\bigO(\sqrt{x}),
    \]
    where $\gamma = \lim\limits_{n\to\infty}\left(-\log{n}+\sum\limits_{k=1}^n 1/k\right)\approx 0.577216$ is the Euler-Mascheroni constant.
\end{thm}

Next, let $F(x)=\DirichletD(x)-2\DirichletD(x/2)$. Then $\Fseq_n=F(n)$ and, via Theorem~\ref{thm:dirichlet} we have $F(x)=x\log2+\bigO(\sqrt{x})$. Restricting to integers $n$, we obtain the following result.

\begin{prop}\label{prop:floor_seq_asymp}
    For $n\in\N$, 
    $\Fseq_n 
        = n\log2+\bigO\left(\sqrt{n}\right)$,  
    and hence $\lim\limits_{n\to\infty}\frac{1}{n}\Fseq_n=\log{2}\approx0.693147$.
\end{prop}

\subsection{The ceiling sequence}
As was mentioned in the introduction, for $n\in\N$, $\Cseq_n$ is equal to the number of odd integers in the ceiling sequence 
\begin{equation}\label{seq:ceiling}
    \Ceil{\frac{n}{1}}, \Ceil{\frac{n}{2}}, \dots, \Ceil{\frac{n}{n}},
\end{equation}
where the ceiling function is defined, for $x\in\R$, by $\Ceil{x}=\min\{z\in\Z : z\ge x\}$. In this subsection, we will find a relation between $\Cseq_n$ and $\Fseq_{n-1}$ which, along with results from the previous subsection, will lead to exact and asymptotic formulas for $\Cseq_n$. 

\begin{prop}\label{prop:ceiling_seq}
    For $n\in\N$,
    \[
        \Cseq_n
        = n-\Fseq_{n-1} = \sum\limits_{d=2}^n \Floor{\frac{n}{d}}(-1)^d.
    \]
\end{prop}
\begin{proof}
We begin by showing that $\Ceil{n/k}=\Floor{(n-1)/k}+1$ for all $k\in\N$. As in Lemma~\ref{lem:consec_flr}, via the division algorithm, there are $q,r\in \Z$ for which $n=kq+r$ with $0\leq r< k$. 

If $r=0$, then $n/k=q$ and so $\Ceil{n/k}=\Ceil{q}=q$. We also have
\[
    \flrf{n-1}{k}=\left\lfloor\frac{n}{k}-\frac{1}{k}\right\rfloor
    =\left\lfloor q-\frac{1}{k}\right\rfloor
    =q-1.
    \]
Thus, $\Ceil{n/k}=\Floor{(n-1)/k}+1$.

If $r\ne0$, then $1\le r<k$. We have 
\[
    \clf{n}{k}
    =\clf{kq+r}{k}
    =\left\lceil q+\frac{r}{k}\right\rceil 
    = q+1.
\] 
and
\[
    \flrf{n-1}{k}
    =\flrf{kq+r-1}{k}
    =\left\lfloor q+\frac{r-1}{k}\right\rfloor
    =q,
\] 
with the final equality following from the fact that $1\le r<k$. Thus, $\Ceil{n/k}=\Floor{(n-1)/k}+1$.

Since we have $\Ceil{n/k}=\Floor{(n-1)/k}+1$ for all $k\in\N$, for each pair of integers $\left(\Ceil{n/k},\Floor{(n-1)/k}\right)$, exactly one integer is odd. Thus, the total number of odd integers in the two sequences 
\[
    \Ceil{\frac{n}{1}},\Ceil{\frac{n}{2}},\dots,\Ceil{\frac{n}{n}}
    \text{ and }
    \Floor{\frac{n-1}{1}},\Floor{\frac{n-1}{2}},\dots,\Floor{\frac{n-1}{n}}
\] 
is $n$. The number of odd integers in the first sequence is $\Cseq_n$. Since the last term in the second sequence is $\Floor{(n-1)/n}=0$, the number of odd integers in the second sequence is $\Fseq_{n-1}$. This means $\Fseq_{n-1}+\Cseq_n=n$. The stated result the follows from Proposition~\ref{prop:floor_seq}.
\end{proof} 

We immediately obtain the following asymptotic formula for $\Cseq_n$.

\begin{cor}\label{cor:ceiling_seq_asymp} 
For $n\in\N$, $\Cseq_n
        = (1-\log2)n+\bigO\left(\sqrt{n}\right)$,
    and hence $
        \lim\limits_{n\to\infty}\frac{1}{n}\Cseq_n
        =1-\log2 \approx 0.306853$.
\end{cor}
\begin{proof}
    Combining Proposition~\ref{prop:ceiling_seq} with the asymptotic formula for $\Fseq_n$ in Proposition~\ref{prop:floor_seq_asymp}, we find 
    \[
        \Cseq_n 
        = n - \Fseq_{n-1} 
        = n - (n-1)\log2 + \bigO\left(\sqrt{n-1}\right) 
        = n-n\log2 + \bigO\left(\sqrt{n}\right).
    \]
    The limit result follows.
\end{proof}
\begin{remark}\label{rmk:an-and-bn}
    If we just wanted an asymptotic formula for $\Cseq_n$ without finding an exact formula first, we could have used the asymptotic formula for $\Fseq_n$, the observation that
    \[
        \Ceil{n/k} = 
        \begin{dcases*} 
            \Floor{n/k} & if $k\mid n$ 
            \\ \Floor{n/k}+1 & if $k\nmid n$,
        \end{dcases*}
    \]
    and the following lemma.
    \begin{lem}\label{lem:num_divisors_of_n}
        For $n\in\N$, the number of positive divisors of $n$ is at most $2\sqrt{n}$.
    \end{lem}
    \begin{proof}
        Suppose $n=de$ for positive integers $d\le e$. Then $1\le d\le \sqrt{n}$. It follows that there are at most $\sqrt{n}$ pairs of divisors $d,e$, and thus at most $2\sqrt{n}$ positive divisors of $n$.
    \end{proof}
    
    Thus, whenever $k$ does not divide $n$, exactly one of $\Floor{n/k}$ and $\Ceil{n/k}$ is odd. When $k$ does divide $n$, then either both $\Floor{n/k}$ and $\Ceil{n/k}$ are odd, or neither is. Since there are at most $2\sqrt{n}$ divisors $d$ of $n$, we have 
    \[
        \Fseq_n+\Cseq_n \in [n-2\sqrt{n},n+2\sqrt{n}].
    \]
    Thus $\Fseq_n+\Cseq_n = n + \bigO\left(\sqrt{n}\right)$,
    from which we conclude that $\Cseq_n=n-\Fseq_n+\bigO\left(\sqrt{n}\right) = (1-\log2)n+\bigO\left(\sqrt{n}\right)$. This gives an alternate proof of Corollary~\ref{cor:ceiling_seq_asymp}.
\end{remark}

\section{The nearest integer sequence}\label{sec:Rseq_n}
As was mentioned in the introduction, for $n\in\N$, $\Rseq_n$ counts the number of odd integers in the nearest integer (rounding) sequence 
\begin{equation}\label{seq:nint-seq}
    \nint\left(\frac{n}{1}\right), \nint\left(\frac{n}{2}\right), \dots, \nint\left(\frac{n}{n}\right),
\end{equation}
where the nearest integer function is defined, for $x\in\R$, by 
\[\nint(x)=\max\{z\in\Z : |z-x|\le |z'-x|\text{ for all } z'\in\Z\}.\]
(In other words, we round to the nearest integer, and half-integers are rounded up.) In this subsection, we will find an exact formula for $\Rseq_n$. Then, with that formula and the solution to Gauss's circle problem, we will find an asymptotic formula for $\Rseq_n$.

We should first note that $\nint(n/k)=\Floor{n/k+1/2}$, and thus Sequence~\eqref{seq:nint-seq} is equal to the sequence \begin{equation}\label{seq:rounding-seq}
    \Floor{\frac{n}{1}+\frac{1}{2}}, \Floor{\frac{n}{2}+\frac{1}{2}}, \Floor{\frac{n}{3}+\frac{1}{2}}, \dots, \Floor{\frac{n}{n}+\frac{1}{2}}.
\end{equation}

Just as we did in obtaining an exact formula for $\Fseq_n$ (Proposition~\ref{prop:floor_seq}), we begin with a lemma.
\begin{lem}\label{lem:round_count}
    Fix $n\in\N$. For $k\in\N$, let $g(k)$ equal the number of terms in Sequence~\eqref{seq:rounding-seq} that are greater than or equal to $k$. Then $g(1)=n$ and, for $k\ge2$, $g(k)=\Floor{2n/(2k-1)}$.
\end{lem} 
\begin{proof}
Consider an arbitrary term $\Floor{n/d+1/2}$ in the Sequence~\eqref{seq:rounding-seq}, so that $1\le d\le n$. We first observe that $\Floor{n/d+1/2}\ge\Floor{n/n+1/2} = 1$. Hence, $g(1)=n$.

Next, for $k\ge2$, we have either $d\le 2n/(2k-1)$ or $d>2n/(2k-1)$. We'll consider these cases separately. 

If $d\leq 2n/(2k-1)$, then $k-1/2 \leq n/d$ which implies $k \leq n/d+1/2$. Since $k$ is an integer, we take the floor of both sides to find $k \leq\Floor{n/d+1/2}$.

If $d>2n/(2k-1)$, then $k-1/2>n/d$, which implies $k > n/d+1/2 \geq \Floor{n/d+1/2}$. 

Thus, for $k\ge2$, we have $\Floor{n/d+1/2}\geq k$ for $d \in\left\{1, 2, \dots,\Floor{2n/(2k-1)}\right\}$. Since $\Floor{2n/(2k-1)}\le n$, we conclude that $g(k)=\Floor{2n/(2k-1)}$.
\end{proof}

With this lemma, we can now find an exact formula for $\Rseq_n$.

\begin{prop}\label{prop:rounding_seq}
    For $n\in\N$,
    \[
        \Rseq_n 
        = -n+\sum\limits_{d=1}^{n}\Floor{\frac{2n}{2d-1}}(-1)^{d+1}.
    \]
\end{prop}
\begin{proof}
We wish to count the number of odd terms in Sequence~\eqref{seq:nint-seq}. By Lemma~\ref{lem:round_count}, there are $n$ terms that are at least 1, and, for $k\ge2$, there are $\Floor{2n/(2k-1)}$ terms that are at least $k$. We can use this to count the number of terms that are equal to a given value.

For $k=1$, there are $n-\Floor{2n/3}$ terms that are equal to 1. For $k\ge2$, there are $\Floor{2n/(2k-1)}-\Floor{2n/(2k+1)}$ terms that are equal to $k$. Summing over odd $k$, we find
\begin{align*}
    \Rseq_n 
    &=n-\Floor{\frac{2n}{3}}+\sum\limits_{\substack{3\le k\le n \\ k\text{ odd}}}\Floor{\frac{2n}{2k-1}}-\Floor{\frac{2n}{2k+1}} \\
    &=-n+\sum\limits_{\substack{1\le k\leq n\\ k\text{ odd}}} \flrf{2n}{2k-1}-\flrf{2n}{2k+1}\\
    &=-n+\sum\limits_{d=1}^n\Floor{\frac{2n}{2d-1}}(-1)^{d+1},
\end{align*}
which completes the proof.
\end{proof}

Next, we will work toward an asymptotic formula for the summation in $\Rseq_n$. We will use results of Jacobi and Gauss.

To start, we need some notation. For $n\in\N$, $r\in\Z$, and $m\in\N$, let $\divs_{r,m}(n)$ denote the number of positive divisors of $n$ that are congruent to $r$ modulo $m$. We will be interested in values of this function for $m=4$ and $r$ equal to 1 or 3. For example, with $n=45$, we have $\divs_{1,4}(45)=4$ and $\divs_{3,4}(45)=2$ because the positive divisors of 45 are 1, 3, 5, 9, 15, and 45.

\begin{lem}
    For $n\in\N$,
    \[
        \Rseq_n 
        = -n + \sum\limits_{k=1}^{2n}\left(\divs_{1,4}(k)-\divs_{3,4}(k)\right).
    \]
\end{lem}
\begin{proof}
    Each integer $d$ divides $\Floor{2n/d}$ integers in the interval $[1,2n]$. If we want to count the number of divisors that are 1 modulo 4 for all of the integers from 1 to $2n$, we see that 1 is a divisor of $\Floor{2n/1}$ terms, 5 is a divisor of $\Floor{2n/5}$ terms, 9 is a divisor of $\Floor{2n/9}$ terms, and so on. In other words, we get 
    \[
        \sum\limits_{k=1}^{2n}\divs_{1,4}(k)
        =\Floor{\frac{2n}{1}}+\Floor{\frac{2n}{5}}+\Floor{\frac{2n}{9}}+\dots.
    \] 
    Similarly, if we want to add up the number of divisors that are 3 modulo 4 for all of the integers from 1 to $2n$, we get
    \[
        \sum\limits_{k=1}^{2n}\divs_{3,4}(k)
        =\Floor{\frac{2n}{3}}+\Floor{\frac{2n}{7}}+\Floor{\frac{2n}{11}}+\dots.
    \]
    Observe that the terms in each of the two above summations are eventually all 0, and thus these are finite sums.
    
    Hence,
    \begin{align*}
        \sum\limits_{d=1}^n\Floor{\frac{2n}{2d-1}}(-1)^{d+1} 
        &= \Floor{\frac{2n}{1}}-\Floor{\frac{2n}{3}}+\Floor{\frac{2n}{5}}-\Floor{\frac{2n}{7}}+\cdots+(-1)^{n+1}\Floor{\frac{2n}{2n-1}} \\
        &= \sum\limits_{k=1}^{2n}\left(\divs_{1,4}(k)-\divs_{3,4}(k)\right).
    \end{align*}
    The stated result then follows from the formula for $\Rseq_n$ in Proposition~\ref{prop:rounding_seq}.
\end{proof}

Next, let $\Jacr_2(n)$ be the number of representations of $n$ as a sum of two integer squares. More precisely, $\Jacr_2(n)=\#\{(a,b)\in\Z^2 : a^2+b^2=n\}$.  For example, still with $n=45$, $\Jacr_2(45)=8$ because
\[
    45=(\pm3)^2+(\pm6)^2=(\pm6)^2+(\pm3)^2,
\] 
which is a total of 8 combinations. Jacobi's two-square theorem relates $\Jacr_2(n)$ to the number of divisors of $n$ that are 1 modulo 4 and that are 3 modulo 4.

\begin{thm}[Jacobi's two-square theorem]\label{thm:jacobi}
    For $n\in\N$, $\Jacr_2(n)=4(\divs_{1,4}(n)-\divs_{3,4}(n))$.
\end{thm}

Jacobi proved this theorem, along with theorems about the number of representations of $n$ using four squares, using six squares, and using eight squares, in 1829 with the use of elliptic theta functions. (See \cite{Grosswald1985}.) For another approach, one can use the theory of binary quadratic forms. (We will use results about binary quadratic forms in Section~\ref{sec:applications}. See, e.g., \cite{Dickson1958}.)

One may also prove Theorem~\ref{thm:jacobi} in the context of the Gaussian integers, which is the ring 
\[
    \Z[\ii]=\{a+b\ii : a,b\in\Z,\, \ii^2=-1\}.
\]
The norm of the Gaussian integer $a+b\ii$ is $a^2+b^2$. It follows that $\Jacr_2(n)$ is equal to the number of Gaussian integers with norm equal to a given positive integer $n$. With an understanding of what the prime elements of $\Z[\ii]$ are, along with the fact that $\Z[\ii]$ is a unique factorization domain, one may prove Theorem~\ref{thm:jacobi}. (For details, see, e.g., \cite[Theorem 278]{HardyWright1979} or a wonderful video on the 3Blue1Brown channel on YouTube \cite{3b1b-pi-prime-regularities}.)

Visualizing $\Z[\ii]$ as a lattice in the complex plane (with $a+b\ii\in\Z[\ii]$ corresponding to the point $(a,b)$ in the plane), Jacobi's two-square theorem says that the number of lattice points on a circle of radius $\sqrt{n}$ centered at the origin is $4(\divs_{1,4}(n)-\divs_{3,4}(n))$. In general, the point $(a,b)$ is on a circle of radius $\sqrt{a^2+b^2}$ centered at the origin. Thus, adding the numbers of points on circles of radii $\sqrt{1}, \sqrt{2},\dots,\sqrt{2n}$ gives us the total number of lattice points different from the origin in the interior or on the boundary of a circle of radius $\sqrt{2n}$ centered at the origin. Back to our formula for $\Rseq_n$, we now have
\begin{equation}\label{eqn:Rseq_n-circle}
    \Rseq_n 
    = -n + (1/4)\cdot \#\{(a,b)\in\Z^2 : 0<a^2+b^2\le 2n\}.
\end{equation}

Let $\GaussC(x)$ denote the number of lattice points in the interior or on the boundary of a circle of radius $\sqrt{x}$ centered at the origin. Finding an asymptotic formula for $\GaussC(x)$ is known as the \emph{Gauss circle problem}. Equation~\eqref{eqn:Rseq_n-circle} shows us how to calculate $\GaussC(2n)$ in terms of $\Rseq_n$.

\begin{prop}\label{prop:gauss-R_n}
    For $n\in\N$, the number of lattice points within a circle of radius $\sqrt{2n}$ centered at the origin is $4\Rseq_n+4n+1$. That is, $\GaussC(2n)=4\Rseq_n+4n+1$.
\end{prop}

\begin{figure}
    \centering
    \includegraphics[scale=.5]{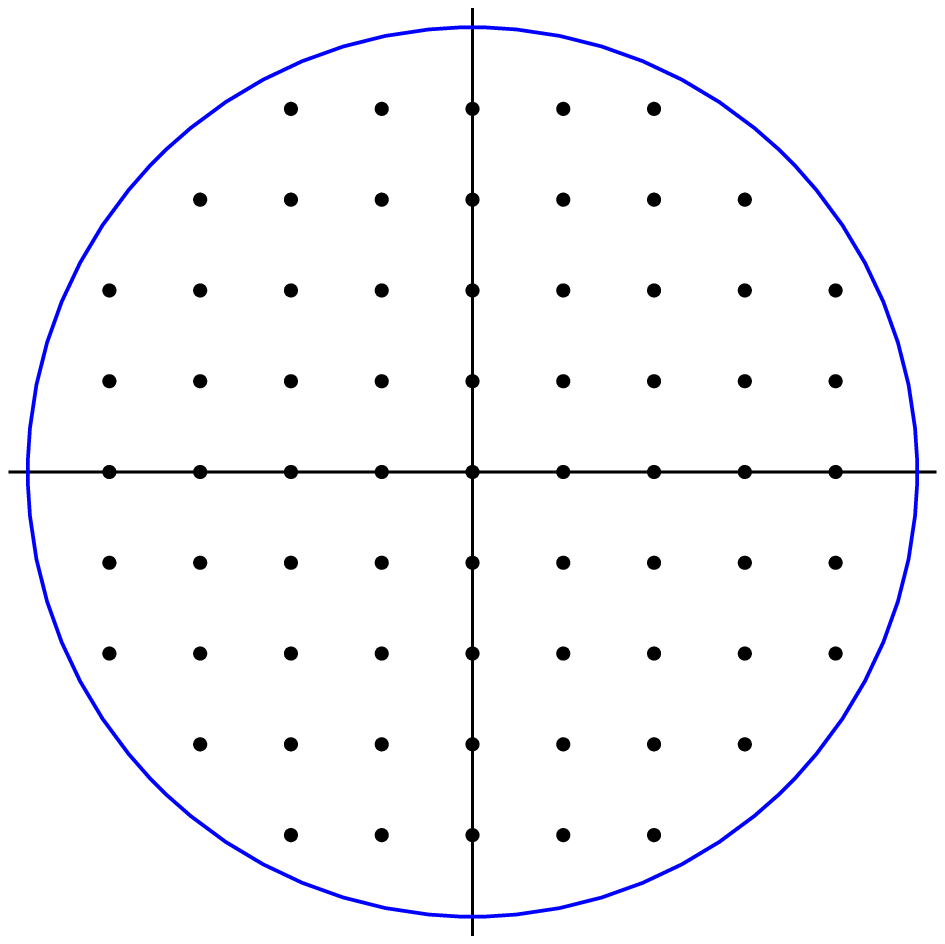}
    \caption{The number of lattice points within and on a circle of radius $6=\sqrt{2\cdot18}$ is, by Proposition~\ref{prop:gauss-R_n},  $\GaussC(2\cdot18)=4\Rseq_{18}+4\cdot18+1=4\cdot10+4\cdot18+1=113$.}
    \label{fig:gauss-circle}
\end{figure}

See Figure~\ref{fig:gauss-circle} for an example with $n=18$, which uses data from Figure~
\ref{fig:first-20-terms}.

Geometrically, we see that $\GaussC(n)$ is a non-decreasing function. For any $n\in\N$, 
\[
    0
    \le \GaussC(2n+2)-\GaussC(2n) 
    = 4\Rseq_{n+1}+4(n+1)+1-4\Rseq_{n}-4n-1 
    = 4\Rseq_{n+1}+4-4\Rseq_n.
\]
Thus, $\Rseq_{n+1}-\Rseq_n\ge-1$. This proves the following corollary.

\begin{cor}
    The sequence $\Rseq_n$ decreases by at most 1 in any step.
\end{cor}
A decrease of 1 from $\Rseq_n$ to $\Rseq_{n+1}$ occurs precisely when $2n+2$ and $2n+1$ each cannot be written as a sum of two integer squares. This occurs when the prime factorizations of $2n+2$ and $2n+1$ each contain some prime which is 3 modulo 4 raised to an odd power. From our data in Figure~\ref{fig:first-20-terms}, we see a decrease by 1 for $n=5$. Observe that $2n+2=12=2^2\cdot3^1$ and $2n+1=11^1$. Neither 11 nor 12 is a sum of two integer squares.

For comparison, there is no bound for the amount in which the sequences $\Fseq_n$ and $\Cseq_n$ can decrease in any step. Indeed, for $k\in\N$ and $n=2^k$, $\Fseq_n-\Fseq_{n-1}=\Cseq_n-\Cseq_{n-1}=-(k-1)$.

We now want an asymptotic formula for $\Rseq_n$. To get there, we use an asymptotic result for $\GaussC(x)$ that is due to Gauss.
\begin{thm}\label{thm:gauss}
    For $x\ge1$, $\GaussC(x) = \pi x + \bigO(\sqrt{x})$.
\end{thm}
This result appears widely in the literature. See, for instance, \cite[Theorem 41]{Rademacher77} or \cite[Chapter 2, Section 7]{Grosswald1985}.

We can now compute an asymptotic formula for $\Rseq_n$.
\begin{prop}\label{prop:rounding_seq_asymp}
    For $n\in\N$, $\Rseq_n = (\pi/2-1)n+\bigO\left(\sqrt{n}\right)$,
    and hence
    $\lim\limits_{n\to\infty}\frac{1}{n}\Rseq_n
        =\pi/2-1
        \approx 0.570796$.
\end{prop}

\begin{proof}
Combining Proposition~\ref{prop:gauss-R_n} and Theorem~\ref{thm:gauss}, we find
\begin{align*}
    \Rseq_n 
    &= -n + (1/4)\cdot \left(\GaussC(2n) - 1\right) \\
    &= -n + (1/4)\cdot \left(2\pi n - 1 + \bigO(\sqrt{2n})\right) \\
    &= -n + (\pi/2) n + \bigO\left(\sqrt{n}\right),
\end{align*}
which proves the first part. The limit behavior follows.
\end{proof}

\begin{remark}
    When we defined $\nint(x)$ in the introduction, we chose to always round up half-integers. Suppose we choose a different convention for rounding half-integers (e.g., always rounding down, or always rounding to the nearest even integer, or something else). Call the new rounding function $\nint'(x)$ and consider the resulting sequence 
    \[
        \Rseq_n'
        =\#\left\{1\le k\le n : \nint'(n/k)\text{ is odd}\right\}.
    \]
    To compute $\left|\Rseq_n'-\Rseq_n\right|$, we need only consider those $k$ for which $n/k$ is a half-integer. But $n/k$ is a half-integer when $n/k=l/2$ for odd $l$, which means $k$ is a divisor of $2n$. By Lemma~\ref{lem:num_divisors_of_n}, $2n$ has at most $2\sqrt{2n}$ divisors. Thus, there are at most $2\sqrt{2n}$ such $k$, and so
    \[
        |\Rseq_n'-\Rseq_n|
        \le 2\sqrt{2n},
    \]
    from which we conclude $\Rseq_n'=\Rseq_n+\bigO\left(\sqrt{n}\right)=n(\pi/2-1)+\bigO\left(\sqrt{n}\right)$.
\end{remark}

\section{Counting by congruence class with shifted floors}\label{sec:generalized}
We can now generalize our results from the preceding sections in three ways.

First, we can write $\nint(x)$ in terms of the floor function, by $\nint(x)=\Floor{x+1/2}$. We may therefore think of $\nint(x)$ as a \emph{$1/2$-shifted floor function} and the corresponding sequence
\[
    \nint\left(\frac{n}{1}\right), \nint\left(\frac{n}{2}\right), \dots, \nint\left(\frac{n}{n}\right) 
    = 
    \Floor{\frac{n}{1}+\frac{1}{2}}, \Floor{\frac{n}{2}+\frac{1}{2}}, \dots, \Floor{\frac{n}{n}+\frac{1}{2}}
\]
as a \emph{$1/2$-shifted floor sequence of length $n$}. We will consider an arbitrary shift $\alpha\in\R$.

Second, we now have sequences of length $n$ where the $k$th term (for $k\in\{1,2,\dots,n\}$) is $\Floor{n/k + \alpha}$. We can offset the numerator $n$ by some real number $\nu$, resulting in a general term of the form $\Floor{(n-\nu)/k+\alpha}$.

Third, in our earlier work we focused on the number of odd terms in each sequence, which means counting the number of terms in each sequence which are congruent to 1 modulo 2. We will instead count the number of terms in each sequence that are congruent to $r$ modulo $m$ for integers $r$ and $m$ with $m\ge1$.

Let's set some notation. For $\alpha,\nu\in\R$, the \emph{$\alpha$-shifted, $\nu$-offset, floor sequence of length $n$} is 
\begin{equation}\label{seq:alpha-shift}
    \Floor{\frac{n-\nu}{1}+\alpha}, \Floor{\frac{n-\nu}{2}+\alpha}, \dots, \Floor{\frac{n-\nu}{n}+\alpha}.
\end{equation}
For $r\in\Z$ and $m\in\N$, let $\Num_{n,\alpha,\nu,r,m}$ equal the number of integers in Sequence~\eqref{seq:alpha-shift} that are congruent to $r$ modulo $m$. In other words, let
\[
    \Num_{n,\alpha,\nu,r,m}
    =\#\left\{1\le k\le n : \Floor{\frac{n-\nu}{k}+\alpha}\equiv r\pmod*{m} \right\}.
\]
Connecting to our earlier work, we see that for $n\in\N$, $\Fseq_n=\Num_{n,0,0,1,2}$ and $\Rseq_n=\Num_{n,1/2,0,1,2}$.

Since every integer is in exactly one congruence class modulo $m$, we immediately see that 
\begin{equation}\label{eqn:sum-is-n}
    \sum\limits_{r=1}^m \Num_{n,\alpha,\nu,r,m}
    =n
\end{equation}
for all $\alpha$, $\nu$, and $m$. If we let $m=1$, we find $\Num_{n,\alpha,\nu,r,1}=n$ for all $\alpha$, $\nu$, and $r$. In what follows, we will assume $m\ge2$. Furthermore, noting that $\Num_{n,\alpha,\nu,r,m} = \Num_{n,\alpha-1,\nu,r-1,m,}$, we may suppose $\alpha\in[0,1)$. Next, since the difference $\Num_{n,\alpha,\nu,r,m}-\Num_{n-1,\alpha,\nu-1,r,m,}$ is 1 or 0 depending on whether $\Floor{(n-\nu)/n+\alpha}$ is congruent to $r$ modulo $m$ or not, we may suppose $\nu\in[0,1)$. Finally, it will be useful to take $r\in[1,m]$, the least positive integer in a given congruence class modulo $m$.

\subsection{A sum of differences of floors}
Our first task is to write $\Num_{n,\alpha,\nu,r,m}$ as a sum of differences of floors. We will generalize the approaches taken in writing $\Fseq_n$ (Proposition~\ref{prop:floor_seq}) and $\Rseq_n$ (Proposition~\ref{prop:rounding_seq}) as summations involving differences of floors.

The following lemma generalizes Lemma~\ref{lem:consec_flr} and Lemma~\ref{lem:round_count}.

\begin{lem}\label{lem:consec-general}
    For $\alpha,\nu\in[0,1)$, $k\in\N$, and $n\in\N$ with $n\alpha\ge\nu$, let $g(k)$ equal the number of integers in an $\alpha$-shifted, $\nu$-offset, floor sequence of length $n$ (Sequence~\eqref{seq:alpha-shift}) that are greater than or equal to $k$. Then
    \[
        g(k) = 
        \begin{dcases*}
            n & if $k=1$,\\
            \Floor{(n-\nu)/(k-\alpha)} & if $k\ge2$.
        \end{dcases*}
    \]
\end{lem}
\begin{proof}
    To start, we have $\nu<1\le n$ for all $n$ Thus $n-\nu>0$, which implies $(n-\nu)/1+\alpha, (n-\nu)/2+\alpha, \dots, (n-\nu)/n+\alpha$ is a decreasing sequence. Looking at the final term, since $\alpha\ge\nu/n$, we have 
    \[
        \frac{n-\nu}{n}+\alpha
        \ge \frac{n-\nu}{n}+\frac{\nu}{n}
        =1.
    \]
    Thus, every term in this sequence is at least 1, and hence their floors are at least 1. This means $g(1)=n$.
    
    Now, suppose $k\ge2$ and let $t=(n-\nu)/(k-\alpha)$. Observe that $0<t<n-\nu\le n$. We will show that the first $\Floor{t}$ terms of Sequence~ \eqref{seq:alpha-shift} are at least $k$ and that the remaining terms are less than $k$. Since Sequence~\eqref{seq:alpha-shift} is a non-increasing sequence, it will suffice to show that
    \begin{equation}\label{eq:ineq-chain}
        \Floor{\frac{n-\nu}{t}+\alpha}\ge k 
        > 
        \Floor{\frac{n-\nu}{t+1}+\alpha}.
    \end{equation}
    
    For the first inequality in Equation~\eqref{eq:ineq-chain}, observe that
    \[
        \frac{n-\nu}{t}+\alpha
        = 
        \frac{n-\nu}{(n-\nu)/(k-\alpha)}+\alpha
        = k.
    \]
    Since $k$ is an integer, we have $\Floor{(n-\nu)/t+\alpha}=\Floor{k}=k\ge k$, as desired.
    
    For the second inequality Equation~\eqref{eq:ineq-chain}, we have
    \[
        \Floor{\frac{n-\nu}{t+1}+\alpha} 
        \le 
        \frac{n-\nu}{t+1}+\alpha 
        < 
        \frac{n-\nu}{t}+\alpha
        =
        k.
    \]
    This completes the proof.
\end{proof}
Note that in the above proof, we considered the cases $k=1$ and $k\ge2$ separately. Our argument for $k\ge2$ doesn't work for $k=1$ because, for $n\alpha>\nu$, our value of $t$ would be $t=(n-\nu)/(1-\alpha)>n$, and we cannot have more than $n$ terms in a sequence of $n$ terms.  This explains why we had to get rid of ``extra'' terms in the summation formula for $\Rseq_n$ (Proposition~\ref{prop:rounding_seq}), which involves $k=1$, $\alpha=1/2$, and $\nu=0$, whereas we had no such adjustment in the summation formula for $\Fseq_n$ (Proposition~\ref{prop:floor_seq}), which has $\alpha=\nu=0$.

In what follows, we will count the number of terms in Sequence~\eqref{seq:alpha-shift} that are congruent to $r$ modulo $m$. Since we have slightly different results for $k=1$ and for $k\ge2$ in Lemma~\ref{lem:consec-general}, we will have slightly different results for congruence classes with $r=1$ and with $2\le r\le m$ in the proposition (and subsequent results) below.

\begin{prop}\label{prop:generalized-sum-of-floors}
    Suppose $m\ge2$ and $\alpha,\nu\in[0,1)$. Then for all $n\in\N$ with $n\alpha\ge\nu$,
    \[
        \Num_{n,\alpha,\nu,1,m}
        = n-\Floor{\frac{(n-\nu)}{1-\alpha}}+\sum\limits_{i\ge0} \left(\Floor{\frac{(n-\nu)}{1+im-\alpha}} - \Floor{\frac{(n-\nu)}{2+im-\alpha}}\right)
    \]
    and, for $2\le r\le m$, 
    \[
        \Num_{n,\alpha,\nu,r,m} 
        = \sum\limits_{i\ge0} \left(\Floor{\frac{(n-\nu)}{r+im-\alpha}} - \Floor{\frac{(n-\nu)}{r+1+im-\alpha}}\right).
    \]
\end{prop}
\begin{proof}
    By Lemma~\ref{lem:consec-general}, we have a formula for $g(d)$, the number of integers in Sequence~\eqref{seq:alpha-shift} that are greater than or equal to a given value $d$. Now, let $G(d)$ equal the number of integers in Sequence~\eqref{seq:alpha-shift} that are equal to a given value $d$. Then $G(d) = g(d)-g(d+1)$.
    
    We first compute $G(1)$ by $G(1) = g(1)-g(2) = n - \Floor{(n-\nu)/2-\alpha}$. Then, for $d\ge2$, 
    \[
        G(d) 
        = g(d)-g(d+1) 
        = \Floor{(n-\nu)/(d-\alpha)}-\Floor{(n-\nu)/(d+1-\alpha)}.
    \]
    
    To compute $\Num_{n,\alpha,\nu,r,m}$, we need to count the number of integers in Sequence~\eqref{seq:alpha-shift} that are congruent to $r$ modulo $m$. Since $1\le r\le m$, we have $\Num_{n,\alpha,\nu,r,m} = G(r) + G(r+m) + G(r+2m) + \dots$.
    
    For $2\le r\le m$, we have
    \[
        \Num_{n,\alpha,\nu,r,m} 
        = \sum\limits_{i\ge0}G(r+im) 
        = \sum\limits_{i\ge0}\left(\Floor{\frac{(n-\nu)}{r+im-\alpha}} - \Floor{\frac{(n-\nu)}{r+1+im-\alpha}}\right).
    \]
    
    For $r=1$, we have
    \[
        \Num_{n,\alpha,\nu,1,m} 
        = \sum\limits_{i\ge0}G(1+im) 
        = n - \Floor{\frac{(n-\nu)}{2-\alpha}} + \sum\limits_{i\ge1}\left(\Floor{\frac{(n-\nu)}{1+im-\alpha}} - \Floor{\frac{(n-\nu)}{2+im-\alpha}}\right).
    \]
    If we include the $i=0$ term in the summation, to make it more closely resemble the formula for $r\ne 1$, we get the stated result.
\end{proof}

\subsection{An asymptotic formula via summation}
Our next task is to evaluate the sum 
\[
    \sum\limits_{i\ge0}\left(\Floor{\frac{(n-\nu)}{r+im-\alpha}}-\Floor{\frac{(n-\nu)}{r+1+im-\alpha}}\right)
\] 
that appears in Proposition~\ref{prop:generalized-sum-of-floors}. 

We'll start by simplifying $\frac{1}{n-\nu}$ times this sum. For $x\in\R$, let $\{x\}$ denote the fractional part of $x$. (I.e., let $\{x\}=x-\Floor{x}$.) Then
\[
    \frac{1}{n-\nu}\sum\limits_{i\ge0}\left(\Floor{\frac{n-\nu}{r+im-\alpha}}-\Floor{\frac{n-\nu}{r+1+im-\alpha}}\right)
    =A_{\alpha,r,m}+\frac{1}{n-\nu}B_{n,\alpha,\nu,r,m}
\]
for the quantities
\[
    A_{\alpha,r,m}=\sum\limits_{i\ge0}\left(\frac{1}{r+im-\alpha}-\frac{1}{r+1+im-\alpha}\right)
\]
and
\[
    B_{n,\alpha,\nu,r,m}=\sum\limits_{i\ge0}\left(\fpf{n-\nu}{r+im-\alpha}-\fpf{n-\nu}{r+1+im-\alpha}\right).
\]

In general, $A_{\alpha,r,m}$ is an alternating series in which the absolute values of the terms decrease to zero. By the Alternating Series Test, $A_{\alpha,r,m}$ converges. $B_{n,\alpha,\nu,r,m}$ is also an alternating series. If its terms, in absolute value, were decreasing, then we would have $B_{n,\alpha,\nu,r,m}=\bigO(1)$ and thus $(1/n) B_{n,\alpha,\nu,r,m}=\bigO(1/n)$. Unfortunately, this isn't the case. 

If we can get an asymptotic formula for $B_{n,\alpha,\nu,r,m}$, then we will have an asymptotic formula for $\Num_{n,\alpha,\nu,r,m}$. To start, we will revisit our results for $\Fseq_n$ and $\Rseq_n$ from earlier in this paper.

\begin{example}\label{ex:floor}
    The floor sequence $\Fseq_n$. 
    For $\alpha=\nu=0$, $r=1$, and $m=2$, we have $\Num_{n,0,0,1,2}=\Fseq_n$ for all $n\in\N$. Then 
\[
    A_{0,1,2}
    =\frac{1}{1}-\frac{1}{2}+\frac{1}{3}-\frac{1}{4}+\dots
    =\sum\limits_{k=1}^\infty (-1)^{k+1}\frac{1}{k}
\]
and
\[
    B_{n,0,0,1,2} 
    =\fpf{n}{1}-\fpf{n}{2}+\fpf{n}{3}-\fpf{n}{4}+\dots
    =\sum\limits_{k=1}^\infty (-1)^{k+1}\fpf{n}{k}.
\]
We see that $A_{0,1,2}=\log2$. (This is the Maclaurin series for $\ln(1+x)$, which converges for $-1<x\le 1$, evaluated at $x=1$.) By Proposition~\ref{prop:floor_seq_asymp}, $\Num_{n,0,0,1,2}=\Fseq_n=n\log2+\bigO\left(\sqrt{n}\right)$. By Proposition~\ref{prop:generalized-sum-of-floors}, $B_{n,0,0,1,2}=\bigO\left(\sqrt{n}\right)$.
\end{example} 

\begin{example}\label{ex:rounding}
    The rounding sequence $\Rseq_n$. For $\alpha=1/2$, $\nu=0$, $r=1$, and $m=2$, we have $\Num_{n,1/2,0,1,2}=\Rseq_n$ for all $n\in\N$. Then
\[
    A_{1/2,1,2}
    =\frac{2}{1}-\frac{2}{3}+\frac{2}{5}-\frac{2}{7}+\dots
    =\sum\limits_{k=1}^\infty(-1)^{k+1}\frac{2}{2k+1}
\]
and
\[
    B_{n,1/2,0,1,2} 
    =\fpf{2n}{1}-\fpf{2n}{3}+\fpf{2n}{5}-\fpf{2n}{7}+\dots
    =\sum\limits_{k=1}^\infty(-1)^{k+1}\fpf{2n}{2k+1}.
\]
We see that $A_{1/2,1,2}=\pi/2$. (This is the Maclaurin series for $\arctan(x)$, which converges for $-1\le x\le 1$, evaluated at $x=1$.) By Proposition~\ref{prop:rounding_seq_asymp}, $\Num_{n,1/2,0,1,2}=\Rseq_n=(\pi/2-1)n+\bigO\left(\sqrt{n}\right)$.
Applying Proposition~\ref{prop:generalized-sum-of-floors}, we find $B_{n,1/2,0,1,2}=\bigO\left(\sqrt{n}\right)$.
\end{example} 

As we will see with the following proposition, we have $B_{n,\alpha,\nu,r,m}=\bigO\left(\sqrt{n}\right)$ in general. For the proof, we will use Dirichlet's hyperbola method. Our method is adapted from the approach given in an answer on Mathematics Stack Exchange \cite{MSE-floor-summation} which proved that, for an increasing sequence of positive integers $b_1,b_2,b_3,\dots$,  
\[
    \sum\limits_{k\le n}\Floor{\frac{n}{b_k}}(-1)^k 
    = n\sum\limits_{k\le n}\frac{1}{b_k}(-1)^k + \bigO\left(\sqrt{n}\right).
\]
The terms in our series are not quite in this form, so we will prove a slightly more general result. To do so, we need a generalized notion of the term ``divides'' to work with real numbers.

\begin{defn}[Real-ly divides]
    Let $a\in\R$ and $b\in\Z$. We say $a$ \emph{real-ly divides} $b$ if there is some $d\in\Z$ for which $\Ceil{da}=b$. We denote this by $a\reallydivs b$, and we say $a$ is a \emph{real divisor} of $b$.
\end{defn}

To see that this definition generalizes the usual definition of ``divides,'' let's suppose that we have $a,b\in\Z$ with $a\mid b$. Then there is some $d\in\Z$ for which $da=b$. Hence, $\Ceil{da}=da=b$, and so $a\reallydivs b$.

\begin{remark}
    We will only consider use this definition of \emph{real-ly divides} with positive numbers. However, if one wants to use this with negative numbers as well, it may be beneficial to modify this definition so that it has some symmetry with positive and negative numbers. One could say a real number $a$ real-ly divides an integer $b$ if there is some integer $d$ for which one of the following holds: either $b\ge0$ and $\Ceil{da}=b$; or $b<0$ and $\Floor{da}=b$. With this, one would additionally have $a\reallydivs b$ if and only if $(-a)\reallydivs b$.
\end{remark}

The key property that we need is the following lemma.

\begin{lem}\label{lem:really-divides-summation}
    Let $n\in\N$, $\nu\in[0,1)$, and $a\in\R$ with $a\ge1$. Then 
    \[
        \Floor{\frac{n-\nu}{a}} = \sum\limits_{\substack{1\le d\le n-\nu \\ a\reallydivs d}} 1.
    \]
\end{lem}
\begin{proof}
    To start, we have $n-\nu>0$ and 
    \[
        \Floor{(n-\nu)/a} 
        = \#\left\{ka  : k\in\Z,\, k\ge 1,\, ka\le n-\nu\right\}.
    \]
    Since $a\ge1$, $\Ceil{k_1a}\ne\Ceil{k_2a}$ for all integers $k_1\ne k_2$. Thus,
    \[
        \Floor{(n-\nu)/a} 
        = \#\left\{\Ceil{ka}  : k\in\Z,\, k\ge 1,\, ka\le n-\nu\right\}.
    \]
    But this is just the cardinality of the set of numbers that $a$ real-ly divides. We therefore have
    \[
        \Floor{(n-\nu)/a} 
        = \#\left\{d\in\Z  : 1\le d\le n-\nu,\, a\reallydivs d\right\} = \sum\limits_{\substack{1\le d\le n-\nu \\ a\reallydivs d}} 1,
    \]
    as desired.
\end{proof}

We can now apply Dirichlet's hyperbola method to show $B_{n,\alpha,\nu,r,m}=\bigO\left(\sqrt{n}\right)$.

\begin{prop}\label{prop:sum-of-floors-to-sum-of-fractions}
    For any increasing sequence $b_0,b_1,b_2,\dots$ of positive real numbers with the property that $b_k\ge1$ and $b_k\ge k$ for all $k$, we have
    \[
        \sum\limits_{k\le n-\nu}\Floor{\frac{n-\nu}{b_k}}(-1)^k = n\sum\limits_{k\le n-\nu}\dfrac{(-1)^k}{b_k}+\bigO\left(\sqrt{n}\right).
    \]
\end{prop}
\begin{proof}
    Let $f(n)=\displaystyle\sum\limits_{k\le n-\nu} \Floor{\frac{n-\nu}{b_k}}(-1)^k$. To start, by Lemma~\ref{lem:really-divides-summation}, we have
    \[
        f(n) 
        = \sum\limits_{k\le n-\nu}(-1)^k\sum\limits_{\substack{d\le n-\nu \\ b_k\reallydivs d}} 1.
    \]
    
    Changing the order of summation,
    \[
        f(n)
        = \sum\limits_{d\le n-\nu}\sum\limits_{\substack{b_k\reallydivs d \\ k\le n-\nu}} (-1)^k.
    \]
    We are summing over $d\le n-\nu$ and real divisors $b_k$ of $d$. Thus, $b_k\le n-\nu$. As we have assumed that $b_k\ge k$, this implies $k\le n-\nu$. Thus, we need not explicitly state that $k\le n-\nu$. 
    We have
    \[
        f(n)
        = \sum\limits_{d\le n-\nu}\sum\limits_{b_k\reallydivs d} (-1)^k.
    \]
    If $b_k\reallydivs d$, then $d-1<b_kd'\le d$. Thus, instead of summing over $d\le n-\nu$, we may sum over $d'$ and $k$ such that $b_kd'\le n-\nu$:
    \[
        f(n) 
        = \sum\limits_{b_kd'\le n-\nu}(-1)^k.
    \]

    We will now use Dirichlet's hyperbola method. Consider the region $R$ in the first quadrant of the $xy$-plane that is bounded by the hyperbola $xy=n-\nu$, the line $x=1$, and the line $y=1$. Let $A>0$. We split the region $R$ into 3 subregions: $R_1$, the portion of $R$ which lies above the line $y=(n-\nu)/A$; $R_2$, the portion of $R$ which lies to the right of the line $x=A$; and $R_3$, which is the rectangle $[1,A]\times[1,(n-\nu)/A]$. (See Figure~\ref{fig:hyperbola}.)
    
    \begin{figure}
        \centering
        \hspace{1in} 
        \includegraphics[scale=.7]{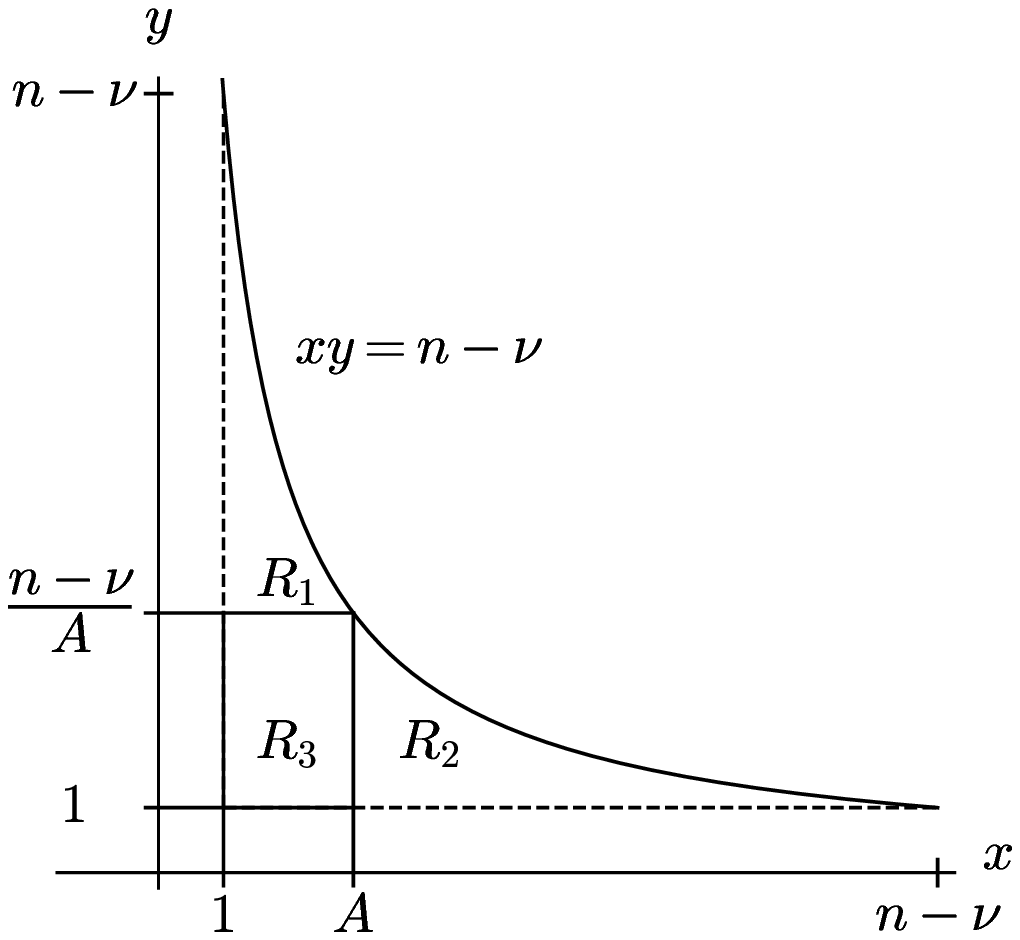}
        \caption{The regions $R_1$, $R_2$, $R_3$}
        \label{fig:hyperbola}
    \end{figure}
    
    It follows that for each combination of $b_k$ and $d'$ such that $b_kd'\le n-\nu$, there is a point $(x,y)=(d',b_k)$ in $R$. Hence, this point is in exactly one of $R_1$, $R_2$, $R_3$.
    
    We can sum over points $(d',b_k)$ in $R_1\cup R_3$ and $R_2\cup R_3$, and then subtract a summation over $R_3$ because we have double counted. We have
    \begin{equation}\label{eqn:dirichlet-eq-5}
        f(n) 
        = 
        \sum\limits_{d'\le A}\sum\limits_{b_k\le (n-\nu)/d'}(-1)^k 
        + \sum\limits_{b_k\le (n-\nu)/A}\sum\limits_{d'\le (n-\nu)/b_k}(-1)^k 
        - \sum\limits_{d'\le A}\sum\limits_{b_k\le (n-\nu)/A} (-1)^k.
    \end{equation}
    Since $b_0,b_1,b_2,\dots$ is an increasing sequence, $\sum\limits_{b_k\le x}(-1)^k\in\{0,1\}$. This summation is $\bigO(1)$. The first and third double sums in Equation~\eqref{eqn:dirichlet-eq-5} are $\bigO(A)$. Hence,
    \[
        f(n) 
        = \sum\limits_{b_k\le (n-\nu)/A}\sum\limits_{d'\le (n-\nu)/b_k}(-1)^k + \bigO(A)
        =\sum\limits_{b_k\le (n-\nu)/A}(-1)^k\Floor{\frac{(n-\nu)}{b_k}}+\bigO(A).
    \]
    Then, since $\Floor{(n-\nu)/b_k}=(n-\nu)/b_k + \bigO(1)$,
    \[
        f(n) =(n-\nu) \sum\limits_{b_k\le (n-\nu)/A} \frac{(-1)^k}{b_k} + \bigO\left(A+\frac{(n-\nu)}{A}\right).
    \]
    Taking $A=\sqrt{(n-\nu)}$, we have
    \[
        f(n) 
        = (n-\nu)\sum\limits_{b_k\le \sqrt{(n-\nu)}}\frac{(-1)^k}{b_k}+\bigO\left(\sqrt{n-\nu}\right)
    \]
    All that remains to do is modify the summation so that we sum over $k\le n-\nu$.  Let $K_n=\min\{k : b_k>\sqrt{n-\nu}\}$. Then
    \[
        \left|\sum\limits_{K_n\le k\le n}\dfrac{(-1)^k}{b_k}\right|
        <\left|\frac{(-1)^{K_n}}{b_{K_n}}\right|
        =\dfrac{1}{b_{K_n}}<\dfrac{1}{\sqrt{n-\nu}}
        =\bigO\left(\dfrac{1}{\sqrt{n-\nu}}\right),
    \]
    where the first inequality holds because we have a finite alternating series.
    Thus, 
    \begin{align*}
        \sum\limits_{b_k\le\sqrt{n-\nu}}\dfrac{(-1)^k}{b_k}
        &= \sum\limits_{k< K_n}\dfrac{(-1)^k}{b_k} \\
        &= \sum\limits_{k\le n-\nu}\dfrac{(-1)^k}{b_k} - \sum\limits_{K_n\le k\le n-\nu}\dfrac{(-1)^k}{b_k} \\
        &= \sum\limits_{k\le n-\nu}\dfrac{(-1)^k}{b_k}+\bigO\left(1/\sqrt{n-\nu}\right).
    \end{align*}
    Therefore,
    \begin{align*}
        f(n)
        &= (n-\nu)\sum\limits_{b_k\le \sqrt{n-\nu}}\frac{(-1)^k}{b_k}+\bigO\left(\sqrt{n-\nu}\right) \\
        &= (n-\nu)\sum\limits_{k\le n-\nu}\dfrac{(-1)^k}{b_k}+\bigO\left(\sqrt{n-\nu}\right) \\
        &= (n-\nu)\sum\limits_{k\le n-\nu}\dfrac{(-1)^k}{b_k}+\bigO\left(\sqrt{n}\right).
    \end{align*}
    Since $\nu$ times the convergent alternating sum is $\bigO(1)$, we get the stated result.
\end{proof}

Combining Proposition~\ref{prop:generalized-sum-of-floors} and Proposition~\ref{prop:sum-of-floors-to-sum-of-fractions}, we obtain an asymptotic formula for $\Num_{n,\alpha,\nu,r,m}$. We'll first give the result for $2\le r\le m$, followed by a slight modification to get the result for $r=1$. (In the proof of Proposition~\ref{prop:generalized-sum-of-floors}, we saw that counting with $r=1$ is slightly different than counting with $r\ne1$. Fortunately, via Equation~\eqref{eqn:sum-is-n}, if we can count for $r=2,\dots,m$, then we get a count for $r=1$ for free.)

\begin{cor}\label{cor:generalized-sum-of-floors}
    For $2\le r\le m$, and $\alpha,\nu\in[0,1)$, and $n\in\N$ with $n\alpha\ge\nu$,
    \[
        \Num_{n,\alpha,\nu,r,m} 
        = n\sum\limits_{i\ge0}\left(\frac{1}{r+im-\alpha}-\frac{1}{r+1+im-\alpha}\right) + \bigO\left(\sqrt{n}\right).
    \]
\end{cor}
\begin{proof}
     For integers $r,m$ with $2\le r\le m$ and for $\alpha\in[0,1)$, define the sequence $b_0,b_1,b_2,\dots$ as follows. For $i\ge0$, let $b_{2i}=(r+im)-\alpha$ and $b_{2i+1}=(r+im)-\alpha+1$. Then, for $\nu\in[0,1)$, 
     \begin{equation}\label{eqn:summation-floor-difference}
        \sum\limits_{k\le n-\nu}\Floor{\frac{n-\nu}{b_k}}(-1)^k
        =
        \sum\limits_{i\ge0}\left(\Floor{\frac{n-\nu}{r+im-\alpha}}-\Floor{\frac{n-\nu}{r+1+im-\alpha}}\right).
    \end{equation}
    (While one series is finite and the other is infinite, the terms in the infinite series are zero for all  $i>(n-\nu-r+\alpha)/m$ and hence there are no convergence issues.)
    
    We wish to show the sequence $b_0,b_1,b_2,\dots$ satisfies the conditions of Proposition~\ref{prop:sum-of-floors-to-sum-of-fractions}. We will show $b_k\ge 1$, $b_k\ge k$, and $b_{k+1}>b_k$ for all $k\ge0$.
    
     If $k$ is even, then $k=2i$ for some $i$ and we have $b_k=(r+km/2)-\alpha$. Since $m\ge2$, we have $b_k\ge k + r-\alpha \ge k$. If $k$ is odd, then $k=2i+1$ for some $i$ and we have $b_k=r+(k-1)m/2-\alpha+1$. Since $m\ge2$, we have $b_k\ge (k-1)+r-\alpha+1\ge k$. Thus, $b_k\ge k$ for all $k\ge0$. Additionally, since $b_0=r-\alpha\ge2-\alpha>1$ and $b_k\ge k\ge1$ for all $k\ge1$, we have $b_k\ge1$ for all $k\ge0$. Finally, for $k$ even, $b_{k+1}-b_k=1$, and for $k$ odd, $b_{k+1}-b_k=m-1\ge1$. Hence, $b_{k+1}>b_k$ for all $k\ge0$.
     
    Since the sequence $b_0,b_1,b_2,\dots$ satisfies the conditions of Proposition~\ref{prop:sum-of-floors-to-sum-of-fractions}, we conclude that 
    \[
        \sum\limits_{k\le n-\nu}\Floor{\frac{n-\nu}{b_k}}(-1)^k 
        = n\sum\limits_{k\le n-\nu}\frac{(-1)^k}{b_k} + \bigO\left(\sqrt{n}\right).
    \]
    Next, since we have an alternating series and $b_k\ge k$ for all $k$, we have 
    \[
        \left|\sum\limits_{k>n-\nu}\frac{(-1)^k}{b_k} \right|
        <\left|\frac{(-1)^{n}}{b_{n}}\right|
        =\frac{1}{b_{n}}
        <\frac{1}{n}
        =\bigO\left(\frac{1}{n}\right).
    \]
    Thus,
    \begin{align*} 
        \sum\limits_{k\le n-\nu}\Floor{\frac{n-\nu}{b_k}}(-1)^k 
        &= n\sum\limits_{k\ge0}\frac{(-1)^k}{b_k}-n\sum\limits_{k>n-\nu}\frac{(-1)^k}{b_k}+\bigO\left(\sqrt{n}\right) \\
        &= n\sum\limits_{k\ge0}\frac{(-1)^k}{b_k} -n\bigO\left(\frac{1}{n}\right)+\bigO\left(\sqrt{n}\right) \\
        &= n\sum\limits_{k\ge0}\frac{(-1)^k}{b_k}+\bigO(1)+ \bigO\left(\sqrt{n}\right).
    \end{align*}
    Combined with Equation~\eqref{eqn:summation-floor-difference}, we conclude
    \[
        \sum\limits_{i\ge0}\left(\Floor{\frac{n-\nu}{r+im-\alpha}}-\Floor{\frac{n-\nu}{r+1+im-\alpha}}\right) 
        = 
        n\sum\limits_{i\ge0}\left(\frac{1}{r+im-\alpha}-\frac{1}{r+1+im-\alpha}\right)+\bigO\left(\sqrt{n}\right).
    \]
    
    Together with Proposition~\ref{prop:generalized-sum-of-floors}, this proves the result for $2\le r\le m$.
\end{proof}

We now use Equation~\eqref{eqn:sum-is-n} and Corollary~\ref{cor:generalized-sum-of-floors} to compute $\Num_{n,\alpha,\nu,r,m}$ for $r=1$.

\begin{cor}\label{cor:generalized-sum-of-floors-r=1}
    For any $m\ge2$, $\alpha,\nu\in[0,1)$, and $n\in\N$ with $n\alpha\ge\nu$,
    \[
        \Num_{n,\alpha,\nu,1,m} 
        = \frac{-\alpha n}{1-\alpha}+n\sum\limits_{i\ge0}\left(\frac{1}{1+im-\alpha}-\frac{1}{2+im-\alpha}\right) + \bigO\left(\sqrt{n}\right).
    \]
\end{cor}
\begin{proof}
    By Equation~\eqref{eqn:sum-is-n}, 
    \[
        \Num_{n,\alpha,\nu,1,m} 
        = n-\sum\limits_{r=2}^m \Num_{n,\alpha,\nu,r,m}.
    \]
    Then, by Corollary~\ref{cor:generalized-sum-of-floors}
    \begin{align*}
        \Num_{n,\alpha,\nu,1,m} 
        & = n - \sum\limits_{r=2}^m \left[n\sum\limits_{i\ge0}\left(\frac{1}{r+im-\alpha}-\frac{1}{r+1+im-\alpha}\right)+\bigO\left(\sqrt{n}\right)\right] \\
        & = n - \sum\limits_{r=2}^m n\sum\limits_{i\ge0}\left(\frac{1}{r+im-\alpha}-\frac{1}{r+1+im-\alpha}\right)+\bigO\left(\sqrt{n}\right) \\
        & = n - n\sum\limits_{i\ge0}\sum\limits_{r=2}^m\left(\frac{1}{r+im-\alpha}-\frac{1}{r+1+im-\alpha}\right) +\bigO\left(\sqrt{n}\right),
    \end{align*}
    where we have changed the order of summation in the last step because we have a finite number of convergent alternating series. We have a telescoping series for each $i\ge0$ which we can manipulate by 
    \begin{align*}
        \Num_{n,\alpha,\nu,1,m} 
        & = n - n\sum\limits_{i\ge0} \left(\frac{1}{2+im-\alpha}-\frac{1}{m+1+im-\alpha}\right)+\bigO\left(\sqrt{n}\right) \\
        & = n -\frac{n}{1-\alpha}+\frac{n}{1-\alpha}- n\sum\limits_{i\ge0} \left(\frac{1}{2+im-\alpha}-\frac{1}{m+1+im-\alpha}\right)+\bigO\left(\sqrt{n}\right) \\
        &= n - \frac{n}{1-\alpha} + n\sum\limits_{i\ge0} \left(\frac{1}{1+im-\alpha}-\frac{1}{2+im-\alpha}\right)+\bigO\left(\sqrt{n}\right).
    \end{align*}
    Note that we have merely inserted $-n/(1-\alpha)+n/(1-\alpha)$ into our expression, and that we have not changed the order of summation in doing so. Hence,
    \[
        \Num_{n,\alpha,\nu,1,m}
        = \frac{-\alpha n}{1-\alpha} + n\sum\limits_{i\ge0} \left(\frac{1}{1+im-\alpha}-\frac{1}{2+im-\alpha}\right)+\bigO\left(\sqrt{n}\right),
    \]
    as desired.
\end{proof}

\subsection{An asymptotic formula via integration}
From integral calculus, we know how to evaluate sums like $1-1/2+1/3-1/4+\dots$ and $1-1/3+1/5-1/7+\dots$ via integration of Maclaurin series as mentioned in Example~\ref{ex:floor} and Example~\ref{ex:rounding}. We will take the same approach to evaluate the summation that appears in Corollary~\ref{cor:generalized-sum-of-floors} and Corollary~\ref{cor:generalized-sum-of-floors-r=1}.
The following proposition shows us how. As with previous work in this section, we will first obtain a result $2\le r\le m$ and then use Equation~\eqref{eqn:sum-is-n} to obtain a result for $r=1$.

\begin{prop}\label{prop:alt-sum-integral-formula}
    For integers $r,m$ with $2\le r\le m$, and for any $\alpha\in [0,1)$,
    \[
        \sum\limits_{i\ge0}\left(\frac{1}{r+im-\alpha}-\frac{1}{r+1+im-\alpha}\right)
        =\int\limits_0^1\dfrac{(1-x)x^{r-1-\alpha}}{1-x^m}\dx.
    \]
\end{prop}
\begin{proof}
    To start, for $\beta=r-1-\alpha$, a positive real number, and for any non-negative integer $k$, let
    \[
        f_k(x)
        =x^{\beta}\sum\limits_{i=0}^k\left(x^{im}-x^{im+1}\right).
    \]
    Then, for all $k\ge0$ we have 
    \[
        |f_k(x)| 
        =\left|x^\beta (1-x)\sum\limits_{i=0}^k x^{im} \right| 
        = \left| x^\beta (1-x) \frac{1-x^{(k+1)m}}{1-x^m} \right|
        = \left| \frac{x^\beta}{1+x+x^2+\dots+x^{m-1}} \right| \left|\left(1-x^{(k+1)m}\right)\right|,
    \]
    a product of absolute values of two functions. Each absolute value is at most 1 for all $x\in[0,1]$. (We're using the fact that $\beta>0$ here.) 
    Thus, $|f_k(x)|\le 1\cdot1=1$ for all $x\in[0,1]$. Since 1 is integrable on $[0,1]$, by the Dominated Convergence Theorem we have
    \[
        \lim\limits_{k\to\infty}\int\limits_0^1 f_k(x)\,\dx 
        = \int\limits_0^1 \lim\limits_{k\to\infty} f_k(x)\,\dx.
    \]
    
    Starting with the integral in the statement of this proposition, and applying this limit and integration interchange, we find
    \begin{align*}
        \int\limits_0^1\dfrac{(1-x)x^{r-1-\alpha}}{1-x^m}\dx
        &= \int\limits_0^1 \lim\limits_{k\to\infty} f_k(x)\dx 
        = \lim\limits_{k\to\infty}\int\limits_0^1 f_k(x)\dx
        = \lim\limits_{k\to\infty} \int\limits_0^1 x^{\beta}\sum\limits_{i=0}^k \left(x^{im}-x^{im+1}\right)\dx \\
        &= \lim\limits_{k\to\infty}\sum\limits_{i=0}^k \left(\frac{x^{\beta+1+im}}{\beta+1+im} - \frac{x^{\beta+2+im}}{\beta+2+im}\right)\Bigg|_{x=0}^{x=1} \\
        & = \lim\limits_{k\to\infty}\sum\limits_{i=0}^k\left(\frac{1}{\beta+1+im}-\frac{1}{\beta+2+im}\right) \\
        & = \sum\limits_{i=0}^\infty\left(\frac{1}{r+im-\alpha}-\frac{1}{r+1+im-\alpha}\right),
    \end{align*}
    as desired.
\end{proof}

Combining Corollary~\ref{cor:generalized-sum-of-floors} and Proposition~\ref{prop:alt-sum-integral-formula}, we obtain the following asymptotic formula for $\Num_{n,\alpha,\nu,r,m}$ with $2\le r\le m$, which we can extend to $r=1$ via Equation~\eqref{eqn:sum-is-n}. This is our main result.

\begin{thm}\label{thm:N_m-asymp}
    Suppose $m\ge1$, $\alpha,\nu\in[0,1)$, and $n\in\N$ with $n\alpha\ge\nu$. Then
    \[
        \Num_{n,\alpha,\nu,1,m} = \frac{-\alpha n}{1-\alpha}+n\int\limits_0^1\dfrac{(1-x)x^{-\alpha}}{1-x^m}\dx + \bigO\left(\sqrt{n}\right),
    \]
    and, for $2\le r\le m$,
    \[
        \Num_{n,\alpha,\nu,r,m} = n\int\limits_0^1\dfrac{(1-x)x^{r-1-\alpha}}{1-x^m}\dx + \bigO\left(\sqrt{n}\right).
    \]
\end{thm}
\begin{proof}
    For $2\le r\le m$, the result follows from Corollary~\ref{cor:generalized-sum-of-floors} and Proposition~\ref{prop:alt-sum-integral-formula}. We need to prove the result for $1\le r\le m$, and $\alpha,\nu\in[0,1)$.
    
    By Equation~\eqref{eqn:sum-is-n},
    \[\Num_{n,\alpha,\nu,1,m} 
        = n - \sum\limits_{r=2}^m\Num_{n,\alpha,\nu,r,m}.\]
    Thus,
    \begin{align*}
        \Num_{n,\alpha,\nu,1,m}
        = n - \sum\limits_{r=2}^m\Num_{n,\alpha,\nu,r,m} 
        &= n - \sum\limits_{r=2}^m \left(n\int\limits_0^1\dfrac{(1-x)x^{r-1-\alpha}}{1-x^m}\dx+\bigO\left(\sqrt{n}\right)\right) \\
        &= n - n\int\limits_0^1\dfrac{x^{-\alpha}(1-x)}{1-x^m}\sum\limits_{r=2}^m x^{r-1}\dx+\bigO\left(\sqrt{n}\right) \\
        &= n - n\int\limits_0^1\dfrac{x^{-\alpha}(1-x)}{1-x^m}\sum\limits_{r=2}^m x^{r-1}\dx+\bigO\left(\sqrt{n}\right).
    \end{align*}
    The finite series inside the integral will cancel with the denominator nicely if we include one more term. We do so as follows:
    \begin{align*}
        \Num_{n,\alpha,\nu,1,m} - n\int\limits_0^1 \dfrac{(1-x)x^{-\alpha}}{1-x^m}\dx &= n - n\int\limits_0^1\dfrac{x^{-\alpha}(1-x)}{1-x^m}\sum\limits_{r=1}^m x^{r-1}\dx+\bigO\left(\sqrt{n}\right) \\
        &=n-n\int\limits_0^1 x^{-\alpha}\dx +\bigO\left(\sqrt{n}\right).
    \end{align*}
    Since $0\le\alpha<1$, this improper integral converges to $1/(1-\alpha)$. Solving for $\Num_{n,\alpha,\nu,1,m}$, we find
    \[
        \Num_{n,\alpha,\nu,1,m} = n - \frac{n}{1-\alpha} + n\int\limits_0^1\dfrac{(1-x)x^{-\alpha}}{1-x^m}\dx + \bigO\left(\sqrt{n}\right),
    \]
    which simplifies to the stated result.
\end{proof}

Thus, $\Num_{n,\alpha,\nu,r,m}$ is asymptotically linear and our formula is independent of $\nu$. We record the corresponding slope below in the following corollary.

\begin{cor}\label{cor:N_m-slope}
    For $\alpha,\nu\in[0,1)$ and $1\le r\le m$,
\[
    \lim\limits_{n\to\infty}\frac{1}{n}\Num_{n,\alpha,\nu,r,m} = 
    \begin{dcases*}
        \frac{-\alpha}{1-\alpha}+\int\limits_0^1\dfrac{(1-x)x^{-\alpha}}{1-x^m}\dx & if $r=1$,\\
        \int\limits_0^1\dfrac{(1-x)x^{r-1-\alpha}}{1-x^m}\dx & if $2\le r\le m$.
    \end{dcases*}
\]
\end{cor}

Via integration, we can compute specific values of $\lim\limits_{n\to\infty}\frac{1}{n}\Num_{n,\alpha,\nu,r,m}$ for $1\le r\le m\le 4$. For $\alpha=\nu=0$, exact and rounded values are in Figure~\ref{fig:shift-0-exact-values} and Figure~\ref{fig:shift-0-decimals}. For $\alpha=1/2$ and $\nu=0$, exact and rounded values are in Figure~\ref{fig:shift-1/2-exact-values} and Figure~\ref{fig:shift-1/2-decimals}.

\begin{figure}
    \centering
    \begin{tabular}{c||c|c|c|c|} 
    \backslashbox{$r$}{$m$} & 1 & 2 & 3 & 4 \\ \hline \hline 
    1 & $1$ & $\log\left(2\right)$ & $\frac{1}{9} \, \sqrt{3} \pi$ & $\frac{1}{8} \, \pi + \frac{1}{4} \, \log\left(2\right)$ \\ \hline 
    2 &  & $-\log\left(2\right) + 1$ & $-\frac{1}{18} \, \sqrt{3} \pi + \frac{1}{2} \, \log\left(3\right)$ & $\frac{1}{8} \, \pi - \frac{1}{4} \, \log\left(2\right)$ \\ \hline 
    3 &  &  & $-\frac{1}{18} \, \sqrt{3} \pi - \frac{1}{2} \, \log\left(3\right) + 1$ & $-\frac{1}{8} \, \pi + \frac{3}{4} \, \log\left(2\right)$ \\ \hline 
    4 &  &  &  & $-\frac{1}{8} \, \pi - \frac{3}{4} \, \log\left(2\right) + 1$ \\ \hline 
    \end{tabular}
    \caption{Values of $\lim\limits_{n\to\infty}\frac{1}{n}\Num_{n,0,0,r,m}$ for $1\le m\le 4$}
    \label{fig:shift-0-exact-values}
\end{figure}

\begin{figure}
    \centering
    \begin{tabular}{c||c|c|c|c|} 
    \backslashbox{$r$}{$m$} & 1 & 2 & 3 & 4 \\ \hline \hline 
    1 & $1.000000$ & $0.693147$ & $0.604600$ & $0.565986$ \\ \hline 
    2 &  & $0.306853$ & $0.247006$ & $0.219412$ \\ \hline 
    3 &  &  & $0.148394$ & $0.127161$ \\ \hline 
    4 &  &  &  & $0.087441$ \\ \hline 
    \end{tabular}
    \caption{$\lim\limits_{n\to\infty}\frac{1}{n}\Num_{n,0,0,r,m}$ for $1\le m\le 4$, rounded to 6 decimal places}
    \label{fig:shift-0-decimals}
\end{figure}

\begin{figure}
    \centering
    \begin{tabular}{c||c|c|c|c|} 
    \backslashbox{$r$}{$m$} & 1 & 2 & 3 & 4 \\ \hline \hline 
    1 & $1$ & $\frac{1}{2} \, \pi - 1$ & $\frac{1}{6} \, \sqrt{3} \pi + \frac{1}{2} \, \log\left(3\right) - 1$ & $\frac{1}{4} \, \pi + \frac{1}{4} \, \sqrt{2} \log\left(3+2\sqrt{2}\right) - 1$ \\ \hline 
    2 &  & $-\frac{1}{2} \, \pi + 2$ & $\frac{1}{6} \, \sqrt{3} {\left(\pi - \sqrt{3} \log\left(3\right)\right)}$ & $\frac{1}{4} \, \pi {\left(\sqrt{2} - 1\right)}$ \\ \hline 
    3 &  &  & $-\frac{1}{3} \, \sqrt{3} \pi + 2$ & $\frac{1}{4} \, \pi - \frac{1}{4} \, \sqrt{2} \log\left(3+2\sqrt{2}\right)$ \\ \hline 
    4 &  &  &  & $-\frac{1}{4} \, \pi {\left(\sqrt{2} + 1\right)} + 2$ \\ \hline 
    \end{tabular}
    \caption{$\lim\limits_{n\to\infty}\frac{1}{n}\Num_{n,1/2,0,r,m}$ for $1\le m\le 4$}
    \label{fig:shift-1/2-exact-values}
\end{figure}

\begin{figure}
    \centering
    \begin{tabular}{c||c|c|c|c|} 
    \backslashbox{$r$}{$m$} & 1 & 2 & 3 & 4 \\ \hline \hline 
    1 & $1.000000$ & $0.570796$ & $0.456206$ & $0.408623$ \\ \hline 
    2 &  & $0.429204$ & $0.357594$ & $0.325323$ \\ \hline 
    3 &  &  & $0.186201$ & $0.162173$ \\ \hline 
    4 &  &  &  & $0.103881$ \\ \hline 
    \end{tabular}
    \caption{$\lim\limits_{n\to\infty}\frac{1}{n}\Num_{n,1/2,0,r,m}$ for $1\le m\le 4$, rounded to 6 decimal places}
    \label{fig:shift-1/2-decimals}
\end{figure}

\section{Applications to finding parity and counting lattice points}\label{sec:applications}
Now that we have a formula for $\Num_{n,\alpha,\nu,r,m}$, we focus on a few applications: computing a floor shift which results in an asymptotic 50/50 split of even and odd terms; and counting lattice points in a few families of ellipses.

\subsection{Shifting for parity}
We return to the case where $m=2$ and consider the problem of determining a shift $\alpha$ so that half of the terms are odd and half are even. It amounts to computing $\alpha$ for which 
\[ 
    \lim\limits_{n\to\infty}\frac{1}{n}\Num_{n,\alpha,\nu,1,2}
    =\lim\limits_{n\to\infty}\frac{1}{n}\Num_{n,\alpha,\nu,2,2}
    =1/2.
\]
By Corollary~\ref{cor:N_m-slope} (with the formula for $r=2$ to avoid the extra term out front), we need $\alpha$ such that
\[
    \int\limits_0^1\dfrac{(1-x)x^{1-\alpha}}{1-x^2}\dx 
    = \int\limits_0^1\dfrac{x^{1-\alpha}}{1+x}\dx 
    = 1/2.
\]
(Note that this is independent of $\nu$.) Since we're using $r=2$, for the remainder of this subsection, we will count even entries instead of odd entries in an $\alpha$-shifted floor sequence.

\begin{figure}
    \centering
    \includegraphics[width=.6\textwidth]{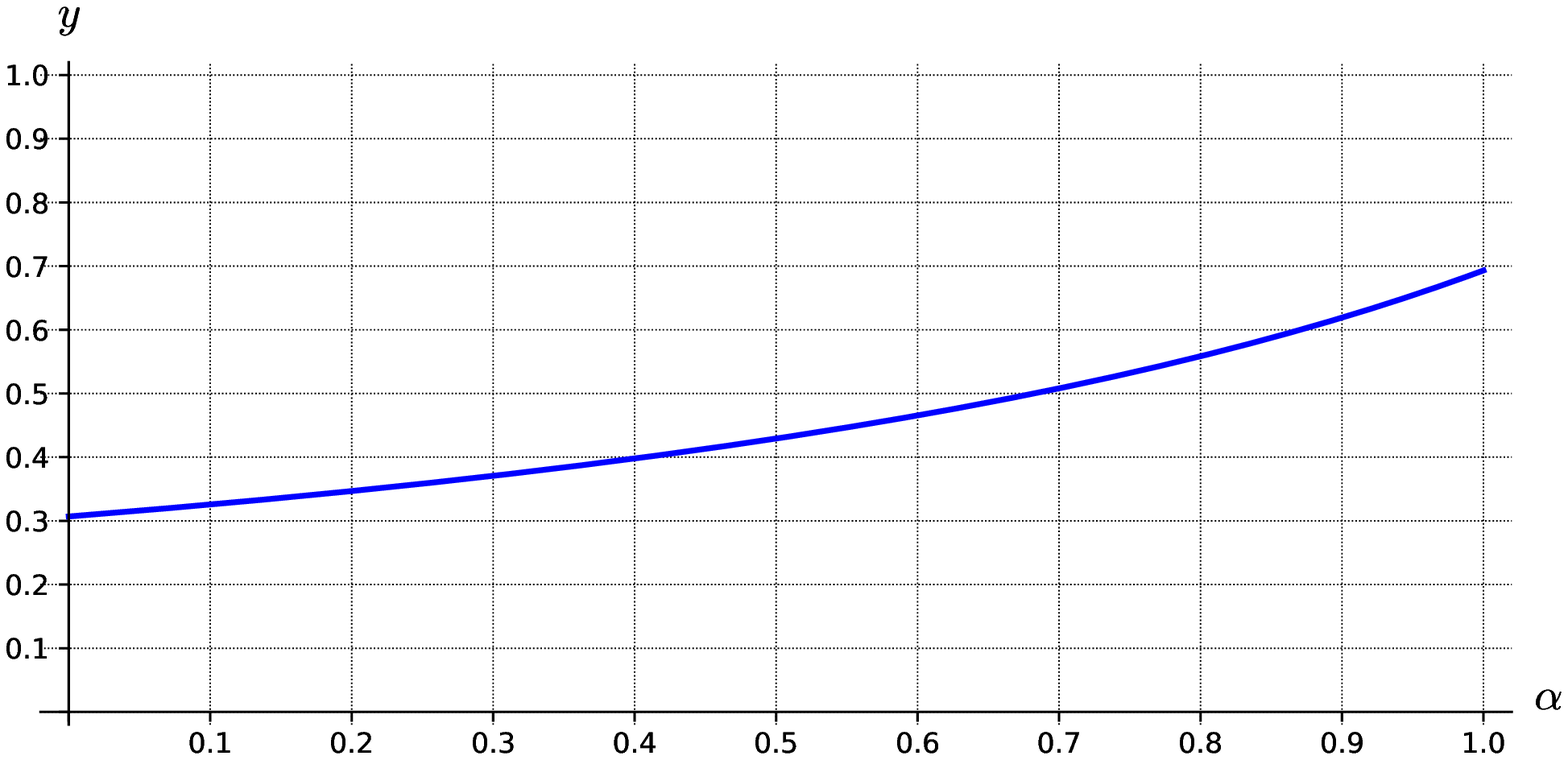}
    \caption{Plot of $\displaystyle y=f(\alpha) = \int\limits_0^1\dfrac{x^{1-\alpha}}{1+x}\dx$ for $\alpha\in[0,1]$. 
    }
    \label{fig:f(x)-plot}
\end{figure}

Let 
\[
    f(\alpha)
    =\int\limits_0^1\dfrac{x^{1-\alpha}}{1+x}\dx.
\]
Then $f(\alpha)$ is the (asymptotic) proportion of terms in an $\alpha$-shifted floor sequence of length $n$ that are even. We immediately see that $f$ is continuous, increasing, and concave up for $\alpha\in[0,1]$.  Furthermore,  $f(0)=1-\log2 < 1/2$ and $f(1)=\log2 > 1/2$. Thus, there is a unique shift $\alpha_0\in(0,1)$ for which $f(\alpha_0)=1/2$. A plot of $f$ appears in Figure~\ref{fig:f(x)-plot}. We see that the value of $\alpha$ for which $f(\alpha)=1/2$ is between $0.6$ and $0.7$.

Computing with Sage \cite{sage}, we can shrink the interval down. We compute
\[
    f(0.682379227335) < 1/2 
    \text{ and } 
    f(0.682379227345) > 1/2.
\]
Hence, we have $f(\alpha_0)=1/2$ for 
\[
    \alpha_0
    \approx 0.68237922734.
\]

\begin{figure}
    \centering
    \includegraphics[width=.6\textwidth]{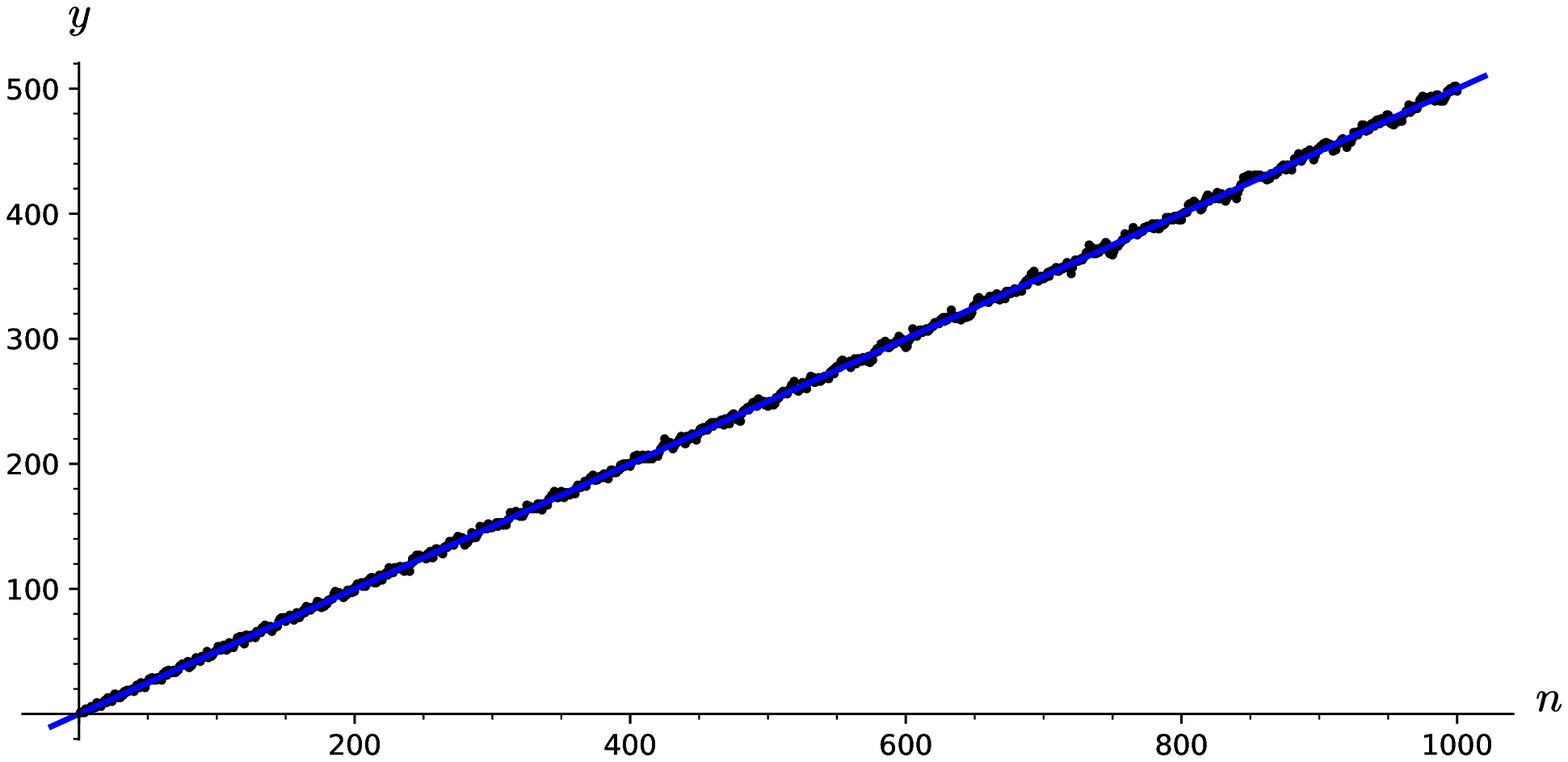}
    \caption{Plot of $y=\Num_{n,\alpha,0,2,2}$ for $\alpha=0.68237933734$ and $1\le n\le 1000$ (black) along with the graph of $y=n/2$ (blue).}
    \label{fig:alpha_0_sequence_plot}
\end{figure}

We can look at some data with this approximate value. A plot of $y=\Num_{n,0.68237922734,0,2,2}$ appears in Figure~\ref{fig:alpha_0_sequence_plot} along with the graph of $y=n/2$. In Figure~\ref{fig:alpha_0_random_table} gives values of $\Num_{n,0.68237922734,0,2,2}$ for random $n$ with $10^5<n<10^6$.

\begin{figure}
    \centering
    \begin{tabular}{|c|c|c|}\hline & & \\
    $n$ & number of even terms 
    & proportion of even terms 
    \\ & & \\
    \hline 
    $ 7015 $ & $ 3503 $ & $49.9359\%$ \\ 
    $ 179220 $ & $ 89632 $ & $50.0123\%$ \\ 
    $ 213788 $ & $ 106901 $ & $50.0033\%$ \\ 
    $ 267093 $ & $ 133562 $ & $50.0058\%$ \\ 
    $ 439675 $ & $ 219839 $ & $50.0003\%$ \\ 
    $ 491213 $ & $ 245600 $ & $49.9987\%$ \\ 
    $ 503521 $ & $ 251741 $ & $49.9961\%$ \\ 
    $ 689325 $ & $ 344631 $ & $49.9954\%$ \\ 
    $ 775294 $ & $ 387629 $ & $49.9977\%$ \\ 
    $ 978029 $ & $ 489010 $ & $49.9995\%$ \\  
    \hline
    \end{tabular}
    \caption{The number and proportion of even terms in an $\alpha$-shifted, $\nu$-offset, floor sequence of length $n$ for $\alpha=0.68237933634$ and various random $n$ with $10^5<n<10^6$.}
    \label{fig:alpha_0_random_table}
\end{figure}

\subsection{Lattice points in selected ellipses, number rings, and plane tilings}
As we saw with Proposition~\ref{prop:gauss-R_n}, the number of lattice points in a circle of radius $\sqrt{2n}$ centered at the origin is $4\Rseq_n+4n+1$, where $\Rseq_n=\Num_{n,1/2,0,1,2}$ is the $1/2$-shifted floor sequence of length $n$. If we think of the plane as being tiled with $1\times1$ squares, then we have a formula for the number of vertices contained in a circle of radius $\sqrt{2n}$ centered at any vertex.

In this subsection, we will find formulas involving $\Num_{n,\alpha,\nu,r,m}$ for the number of lattice points contained in the ellipses
\[
    x^2+y^2=n,\quad x^2+xy+y^2=n,\quad\text{and } x^2+2y^2=n,
\]
for any $n\in\N$. We will also count vertices contained in a circle of radius equal to the square root of any integer for tilings of the plane by squares or triangles.

In Proposition~\ref{prop:generalized-sum-of-floors}, we wrote $\Num_{n,\alpha,\nu,1,m}$ with a summation involving a difference of floors. In what follows, it will be useful to have the formula for $\Num_{n,\alpha,\nu,1,m}$ in the case where $\alpha$ is rational. If we suppose $\alpha=p/q$, then
\begin{equation} \label{eqn:num-rational-alpha}
    \Num_{n,p/q,\nu,1,m}=n - \Floor{\frac{(n-\nu)q}{q-p}} + \sum\limits_{i\ge0}\left(\Floor{\frac{(n-\nu)q}{q+qim-p}}-\Floor{\frac{(n-\nu)q}{2q+qim-p}} \right)
\end{equation}
for all $n\in\N$ with $n\alpha\ge\nu$.

Also, for $n\in\N$, recall the function $\divs_{r,m}(n)$, which counts the number of positive divisors of $n$ that are congruent to $r$ modulo $m$.

\subsubsection{Lattice points in the region $x^2+y^2\le n$}
In Proposition~\ref{prop:gauss-R_n}, we found a formula for the number of lattice points $(x,y)$ contained in the disc $x^2+y^2\le 2n$ in terms of $\Rseq_n$. This is a result for discs with radius equal to the square root of an even number. We'll extend the result to square roots of odd numbers as well, eventually obtaining a formula for the number of lattice points contained in the disc $x^2+y^2\le n$. If we think of the plane as being tiled with $1\times1$ squares, then our formula will give us the total number of vertices that lie in a disc of radius $\sqrt{n}$ centered at one of the vertices.

\begin{prop}\label{prop:circle-radius-n}
    For $n\in\N$, let $F(n)$ be the number of lattice points in a disc of radius $\sqrt{n}$ centered at the origin. Then 
    \[
        F(n) = 
        4\Num_{\Ceil{n/2},1/2,\{n/2\},1,2}+4\Floor{\frac{n}{2}}+1,
    \]
    where $\{n/2\}$ denotes the fractional part of $n/2$.
\end{prop}
\begin{proof}
    Let $F(n)=\#\left\{(x,y)\in\Z^2 : x^2+y^2\le n\right\}$. We'll show that the formula works for even $n$ and for odd $n$.
    
    If $n$ is even, then $n=2k$ for some $k\in\N$.  By Proposition~\ref{prop:gauss-R_n}
    \[
        F(n) 
        = F(2k) 
        = 4\Rseq_k+4k+1=4\Num_{k,1/2,0,1,2}+4k+1=4\Num_{n/2,1/2,0,1,2}+2n+1.
    \]
    Note that $\Ceil{n/2}=k=n/2$, $\{n/2\}=\{k\}=0$, and $4\Floor{n/2}=4\Floor{k}=4k=2n$. This proves the formula for $F(n)$ with $n$ even.
    
    If $n$ is odd, then $n=2k-1$ for some $k\in\N$. By Jacobi's two-square theorem (Theorem~\ref{thm:jacobi}), 
    \[
        F(n) 
        = 1 + 4\sum\limits_{j=1}^n\left(\divs_{1,4}(j)-\divs_{3,4}(j)\right) 
        = 1 + 4\left(\Floor{\frac{n}{1}}-\Floor{\frac{n}{3}}+\Floor{\frac{n}{5}}-\Floor{\frac{n}{7}}+\dots\right).
    \]
    In order to get this alternating floor sum, we use $\alpha=\nu=1/2$, $r=1$, and $m=2$ with Equation~\eqref{eqn:num-rational-alpha}. We have
    \[
        \Num_{k,1/2,1/2,1,2}
        =1-k+\sum\limits_{i\ge0}\left(\Floor{\frac{2k-1}{4i+1}}-\Floor{\frac{2k-1}{4i+3}}\right) 
        = \frac{1-n}{2}+\sum\limits_{i\ge0}\left(\Floor{\frac{n}{4i+1}}-\Floor{\frac{n}{4i+3}}\right)
    \]
    for all $k\ge\nu/\alpha=1$. Then, 
    \[
        F(n) 
        = 1 + 4\left(\Num_{k,1/2,1/2,1,2}+\frac{n-1}{2}\right) 
        = 1 + 4\Num_{(n+1)/2,1/2,1/2,1,2}+2n-2.
    \]
    Note that $\Ceil{n/2}=k=(n+1)/2$, $\{n/2\}=\{k+1/2\}=1/2$, and $4\Floor{n/2}=4\Floor{k-1/2}=4(k-1)=2n-2$. This proves the formula for $F(n)$ with $n$ odd.
\end{proof}

\begin{figure}
    \centering
    \includegraphics[scale=.6]{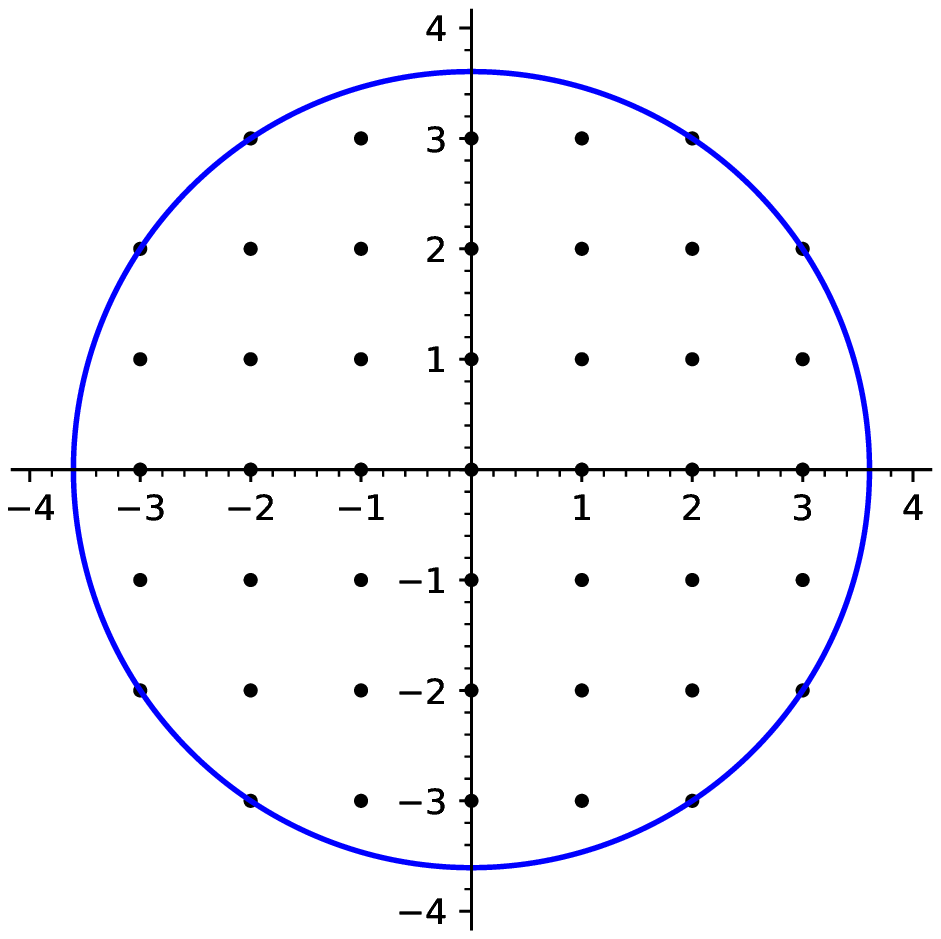}
    \caption{The 45 lattice points in a disc of radius $\sqrt{13}$, which are also the 45 Gaussian integers with norm at most 13.}
    \label{fig:sqrt-13-disc}
\end{figure}

\begin{example}\label{ex:gaussian-ints}
To illustrate Proposition~\ref{prop:circle-radius-n}, we'll compute the number of lattice points in a disc of radius $\sqrt{13}$. Note that $\Ceil{13/2}=7$. The number of lattice points in a disc of radius $\sqrt{13}$ involves the quantity $4\Num_{7,1/2,1/2,1,2}$. To compute this, we examine the length-7 sequence 
\[
    \Floor{\frac{7-1/2}{1}+\frac{1}{2}},\Floor{\frac{7-1/2}{2}+\frac{1}{2}},\dots,\Floor{\frac{7-1/2}{7}+\frac{1}{2}} 
    = 7,3,2,2,1,1,1,
\]
which contains 5 odd terms. Thus, $\Num_{7,1/2,1/2,1,2}=5$. By Proposition~\ref{prop:circle-radius-n}, the number of lattice points in a disc of radius $\sqrt{13}$ is therefore 
\[
    F(13)
    =4\Num_{7,1/2,1/2,1,2}+4\Floor{\frac{13}{2}}+1
    = 4\cdot5+4\cdot6+1=45.
\]
A disc of radius $\sqrt{13}$ appears in Fig.~\ref{fig:sqrt-13-disc}, and one confirms that it contains 45 lattice points.
\end{example}

\begin{remark}
    We mentioned the ring of Gaussian integers, $\Z[\ii]$, earlier. For any $z=a+b\ii\in\Z[\ii]$, the norm of $z$ is $N(z)=N(a+b\ii)=a^2+b^2$. Thus, Proposition~\ref{prop:circle-radius-n} gives a formula for the number of Gaussian integers with norm at most $n$. Following from Example~\ref{ex:gaussian-ints}, we know there are 45 Gaussian integers with norm at most 13. We visualize these Gaussian integers in Figure~\ref{fig:sqrt-13-disc}.
    
    Viewing $\Z[\ii]$ in the complex plane, the Gaussian integers are the vertices for a tiling of the plane by $1\times1$ square tiles. If we draw a circle of radius $\sqrt{n}$, for $n\in\N$, around any lattice point, then Proposition~\ref{prop:circle-radius-n} gives us a formula for the number of vertices contained in the circle. Any other $1\times1$ square tiling of the plane would involve a rotation and/or shift of this tiling. A circle of radius $\sqrt{n}$ centered at any vertex would contain the same number of lattice points. Thus, Proposition~\ref{prop:circle-radius-n} gives us a formula for the number of vertices contained in a circle of radius $\sqrt{n}$, for $n\in\N$, centered at any vertex of any $1\times1$ square tiling of the plane.
\end{remark}

\subsubsection{Lattice points in the region $x^2+xy+y^2\le n$}
Next, we consider the ellipse $x^2+xy+y^2=n$. In general, the quantity $ax^2+bxy+cy^2$, for constants $a,b,c$, is a \emph{binary quadratic form}. We take results about the number of representations of an integer $n$ by a binary quadratic form from \cite{Dickson1958}. 

\begin{prop}[{\cite[Exercise XXII.2]{Dickson1958}}]\label{prop:dickson-x^2+xy+y^2}
    Let $n\in\N$. Then the number of representations of $n=x^2+xy+y^2$, for integers $x,y$, is $6\left(\divs_{1,3}(n)-\divs_{2,3}(n) \right)$.
\end{prop}
\begin{cor}\label{cor:lattice-points-eisenstein-ellipse}
    For $n\in\N$, let $F(n)$ be the number of lattice points contained in the elliptical region $x^2+xy+y^2\le n$. Then
    \[
        F(n)
        = 6\Num_{n,0,0,1,3}+1.
    \]
\end{cor}
\begin{proof}
    Let $f(n)=\#\left\{(x,y)\in\Z^2 : x^2+xy+y^2=n\right\}$. Observe that $f(0)=1$. Thus,
    \begin{equation}\label{eqn:f-step-2}
        F(n) = \sum\limits_{k=0}^n f(k) = 1 + \sum\limits_{k=1}^n f(k).
    \end{equation}
    Next, by Proposition~\ref{prop:dickson-x^2+xy+y^2}, 
    \[
        \sum\limits_{k=1}^n f(k)
        = 6\left(\Floor{\frac{n}{1}}-\Floor{\frac{n}{2}}+\Floor{\frac{n}{4}}-\Floor{\frac{n}{5}}+\Floor{\frac{n}{7}}-\dots\right),
    \]
    which is a finite sum since the floors are eventually zero. Using $\alpha=\nu=0$, $r=1$, and $m=3$ with Equation~\eqref{eqn:num-rational-alpha}, we get
    \[
        \Num_{n,0,0,1,3}
        =\left(\Floor{\frac{n}{1}}-\Floor{\frac{n}{2}}+\Floor{\frac{n}{4}}-\Floor{\frac{n}{5}}+\Floor{\frac{n}{7}}-\dots\right)
    \]
    for all $n\in\N$. We conclude that $F(n) = 1 + 6\Num_{n,0,0,1,3}$,
    as desired.
\end{proof}

\begin{figure}
    \centering
    \includegraphics[scale=.6]{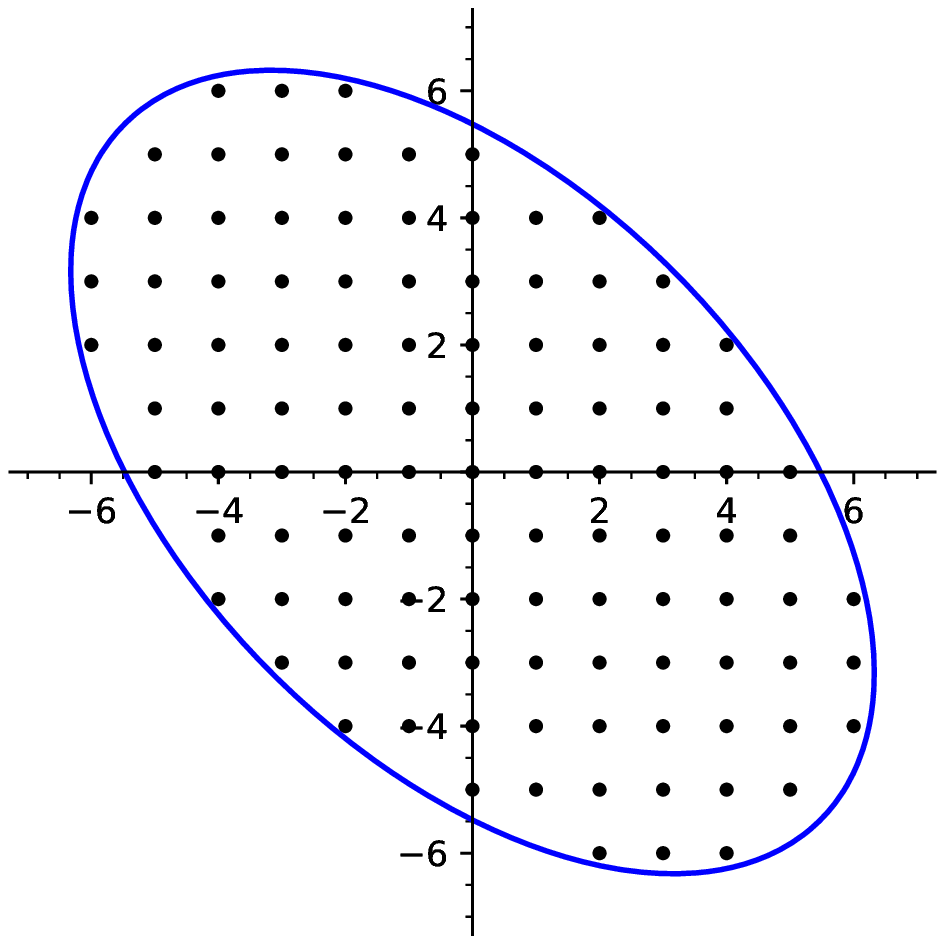}
    \includegraphics[scale=.6]{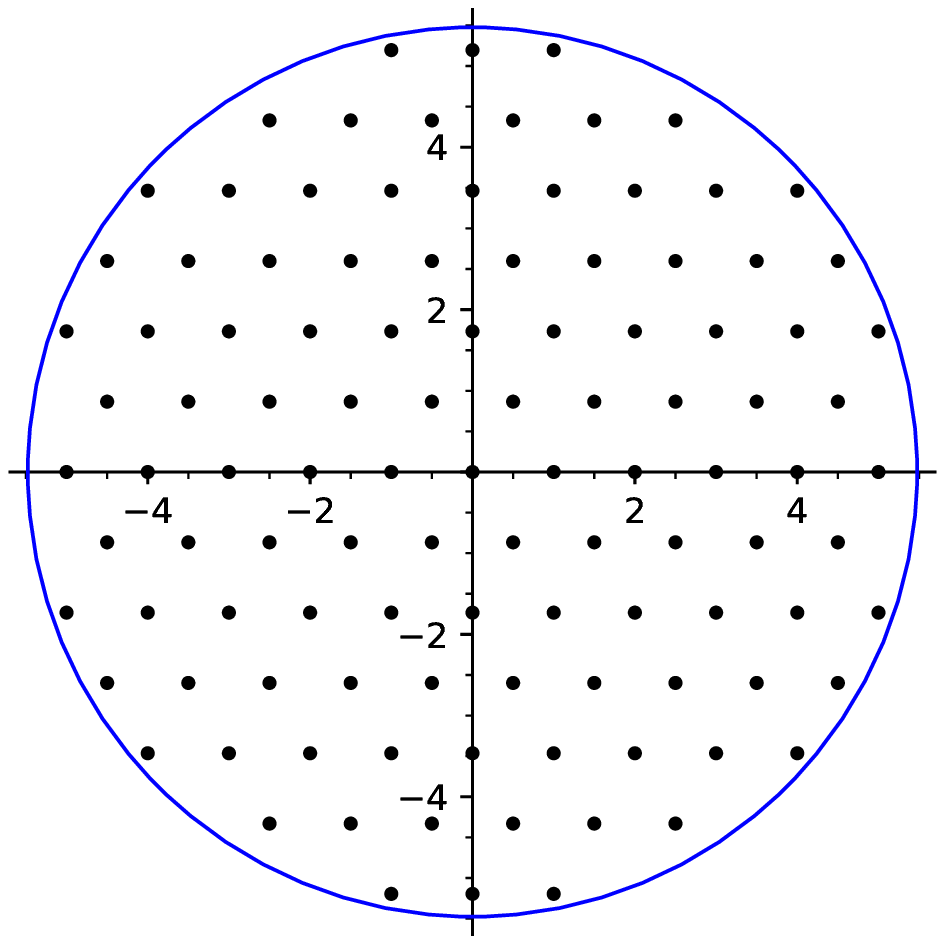}
    \caption{The 109 lattice points within the ellipse $x^2+xy+y^2=30$, and the 109 Eisenstein integers with norm at most 30.}
    \label{fig:ellipse-30}
\end{figure}

\begin{example}\label{ex:ellipse-30}
    We'll compute the number of lattice points contained in the ellipse defined by the equation $x^2+xy+y^2=30$. To do so, we consider the sequence
    \[
        \Floor{\frac{30}{1}},\Floor{\frac{30}{2}},\dots,\Floor{\frac{30}{30}} = 30, 15, 10, 7, 6, 5, 4, 3, 3, 3, 2, 2, 2, 2, 2, 1,1,1,1,1,1,1,1,1,1,1,1,1,1,1.
    \]
    There are 18 terms in this sequence that are congruent to 1 modulo 3. Hence, 
    the ellipse defined by the equation $x^2+xy+y^3=30$ contains $6\Num_{30,0,0,1,3}+1=6\cdot18+1=109$ lattice points. See Figure~\ref{fig:ellipse-30}.
\end{example}

From Theorem~\ref{thm:N_m-asymp} and Figure~\ref{fig:shift-0-exact-values}, we see that $\Num_{n,0,0,1,3}=\frac{1}{9}\sqrt{3}\pi n+\bigO\left(\sqrt{n}\right)$. We immediately obtain the following corollary.

\begin{cor}\label{cor:lattice-points-x2+xy+y2}
    For $n\in\N$, the number of lattice points in the elliptical region $x^2+xy+y^2\le n$ is  $\frac{2}{3}\sqrt{3}\pi n + \bigO\left(\sqrt{n}\right)$.
\end{cor}

\begin{remark} 
    Consider the ring of Eisenstein integers, $\Z[\omega]$, where $\omega=\frac{-1+\sqrt{-3}}{2}$, a primitive 3rd root of unity. For $z=a-b\omega\in\Z[\omega]$, the norm of $z$ is $N(z)=N(a-b\omega)=a^2+ab+b^2$. Thus, Corollary~\ref{cor:lattice-points-eisenstein-ellipse} gives a formula for the number of Eisenstein integers with norm at most $n$. Following from Example~\ref{ex:ellipse-30}, we know there are 109 Eisenstein integers with norm at most 30. We visualize these Eisenstein integers in Figure~\ref{fig:ellipse-30}.
    
    Viewing $\Z[\omega]$ in the complex plane, the Eisenstein integers are the vertices for a plane tiling involving equilateral triangles of side length 1. If we draw a circle of radius $\sqrt{n}$ around any vertex, Corollary~\ref{cor:lattice-points-eisenstein-ellipse} gives us a formula for the number of vertices contained in the circle. Rotating the plane (and hence the tiles) about the center of the circle will leave the number of vertices in the circle unchanged. Thus, Corollary~\ref{cor:lattice-points-eisenstein-ellipse} gives us a formula for the number of vertices contained in a circle of radius $\sqrt{n}$, for $n\in\N$, centered at a vertex of any side-length 1 equilateral triangle tiling of the plane.
\end{remark}

\subsubsection{Lattice points in the region $x^2+2y^2\le n$}
We now consider the ellipse $x^2+2y^2=n$. The result below gives the number of representations of a natural number $n$ in terms of the binary quadratic form $x^2+2y^2$.

\begin{prop}[{\cite[Exercise XXII.1]{Dickson1958}}]
\label{prop:dickson-x^2+2y^2}
    Let $n\in\N$. Then the number of representations of $n=x^2+2y^2$, for integers $x,y$, is $2\left(\divs_{1,8}(n)+\divs_{3,8}(n)-\divs_{5,8}(n)-\divs_{7,8}(n) \right)$.
\end{prop}

\begin{cor}\label{cor:x^2+2y^2}
    Let $F(n)$ be the number of lattice points contained in the elliptical region $x^2+2y^2\le n$. Then $F(0)=1$, $F(1)=3$, $F(2)=5$, $F(5)=11$, and for all $n\in\N$ with $n\ne 1,2,5$, 
    \[
        F(n) 
        = 1+2\Num_{\Ceil{n/4},3/4,\{-n/4\},1,2}+2\Num_{\Ceil{n/4},1/4,\{-n/4\},1,2}+2n+2\Floor{\frac{n}{3}}-4\Ceil{n/4}.
    \]
\end{cor}
\begin{proof}
    Let $f(n) = \#\left\{(x,y)\in\Z^2 : x^2+2y^2=n\right\}$. Observe that $f(0)=1$. To start, we have
    \begin{equation}\label{eqn:f-step-1}
        F(n)
        = \sum\limits_{k=0}^n f(k) 
        = 1 + \sum\limits_{k=1}^n f(k).
    \end{equation}
    Next, by Proposition~\ref{prop:dickson-x^2+2y^2}, 
    \[
        \sum\limits_{k=1}^n f(k) 
        =  2\left(\Floor{\frac{n}{1}}+\Floor{\frac{n}{3}}-\Floor{\frac{n}{5}}-\Floor{\frac{n}{7}}+\Floor{\frac{n}{9}}+\dots\right),
    \]
    which is a finite sum since the floors are eventually zero. We'll need two alternating floor sums here -- one for $\Floor{n/1}-\Floor{n/5}+\Floor{n/9}-\dots$ and one for $\Floor{n/3}-\Floor{n/7}+\Floor{n/11}-\dots$. Each will need $q=4$ and $m=2$. As we did with Proposition~\ref{prop:circle-radius-n}, we will have $\nu\ne0$ here.
    
    For $n\in\N$, let $a=\Ceil{n/4}$ and $b=4a-n$. Then $0\le b<4$. Using $\alpha=p/q=3/4$, $\nu=b/4$, $r=1$, and $m=2$ with Equation~\eqref{eqn:num-rational-alpha}, we get
    \[
        \Num_{a,3/4,b/4,1,2} = a - \Floor{\frac{n}{1}} + \sum\limits_{i\ge0}\left(\Floor{\frac{n}{1+4i}}-\Floor{\frac{n}{5+4i}}\right)
    \]
    for all $a\ge\nu/\alpha=b/3$. 
    Using $\alpha=p/q=1/4$, $\nu=b/4$, $r=1$, and $m=2$ with Equation~\eqref{eqn:num-rational-alpha}, we get
    \[
        \Num_{a,1/4,b/4,1,2} = a - \Floor{\frac{n}{3}} + \sum\limits_{i\ge0}\left(\Floor{\frac{n}{3+4i}}-\Floor{\frac{n}{7+4i}}\right)
    \]
    for all $a\ge\nu/\alpha=b$.
    
    Both equations hold for $a\ge b$. Since $a\ge1$ and $0\le b<4$, this inequality is satisfied for all $a,b$ except for $(a,b)=(1,2), (1,3), (2,3)$, which correspond, respectively, to $n=2$, $n=1$, and $n=5$.
    
    Thus, for $n\in\N$ with $n\ne 1,2,5$, 
    \begin{align*} 
        F(n)
        &=1+2\left(\Num_{a,3/4,b/4,1,2}+\Num_{a,1/4,b/4,1,2}-a+\Floor{\frac{n}{3}}-a+n\right)\\
        &=1+2\Num_{\Ceil{n/4},3/4,\{-n/4\},1,2}+2\Num_{\Ceil{n/4},1/4,\{-n/4\},1,2}+2n+2\Floor{\frac{n}{3}}-4\Ceil{\frac{n}{4}}
    \end{align*} 
    We can fill in the missing values by hand. We compute $f(1)=2$, $f(2)=2$, and $f(5)=0$. And by the formula above with $n=4$,  
    \[
        F(4)
        =1+2\Num_{1,3/4,0,1,2}+2\Num_{1,1/4,0,1,2}+4+2\Floor{4/3}-0
        =1+2\cdot1+2\cdot1+4+2
        =11.
    \]
    Thus, $F(0)=1$, $F(1)=F(0)+f(1)=1+2=3$, $F(2)=F(1)+f(2)=3+2=5$, and $F(5)=F(4)+f(5)=11+0=11$.
\end{proof}
\begin{figure}
    \centering
    \includegraphics[scale=.6]{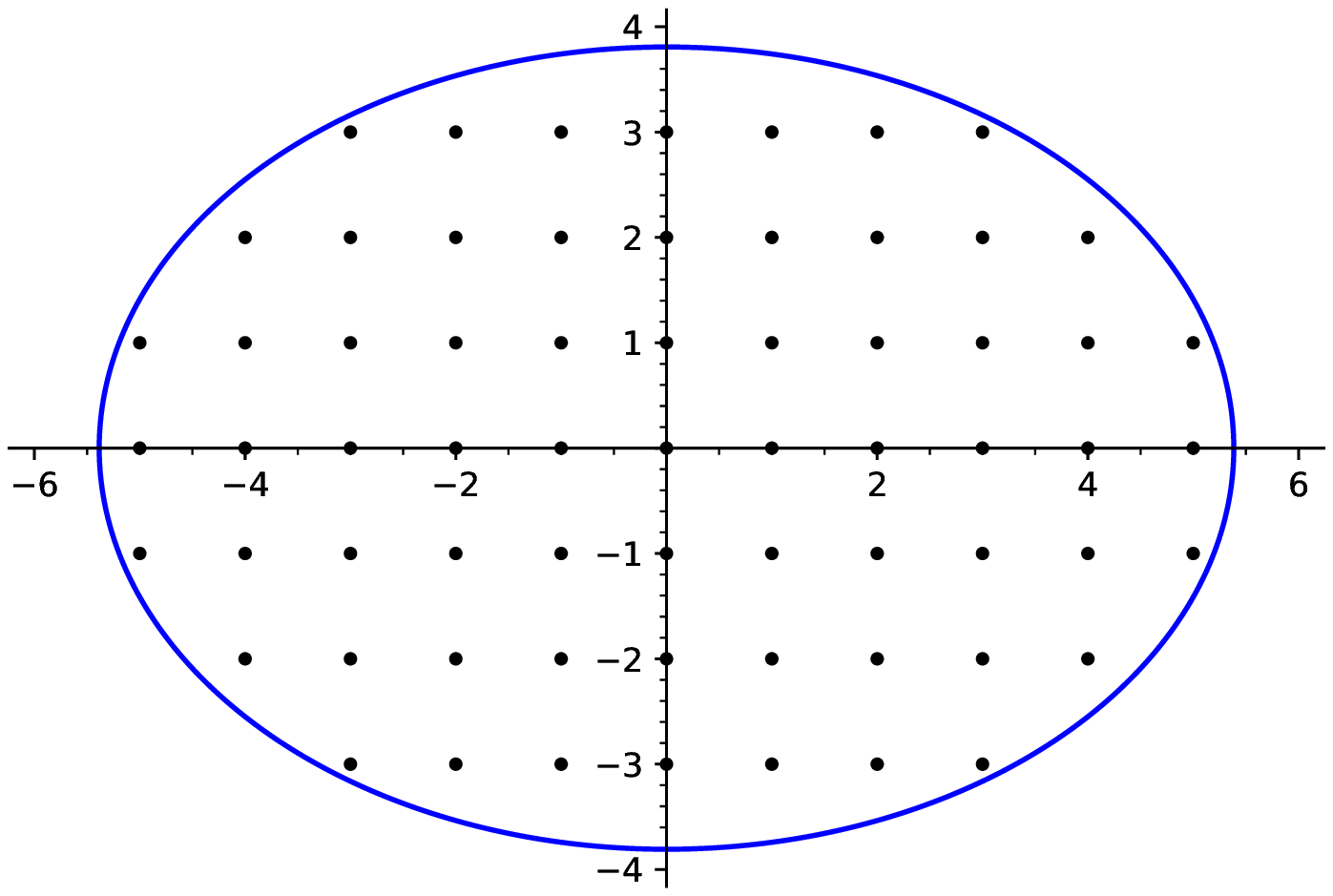}
    \includegraphics[scale=.6]{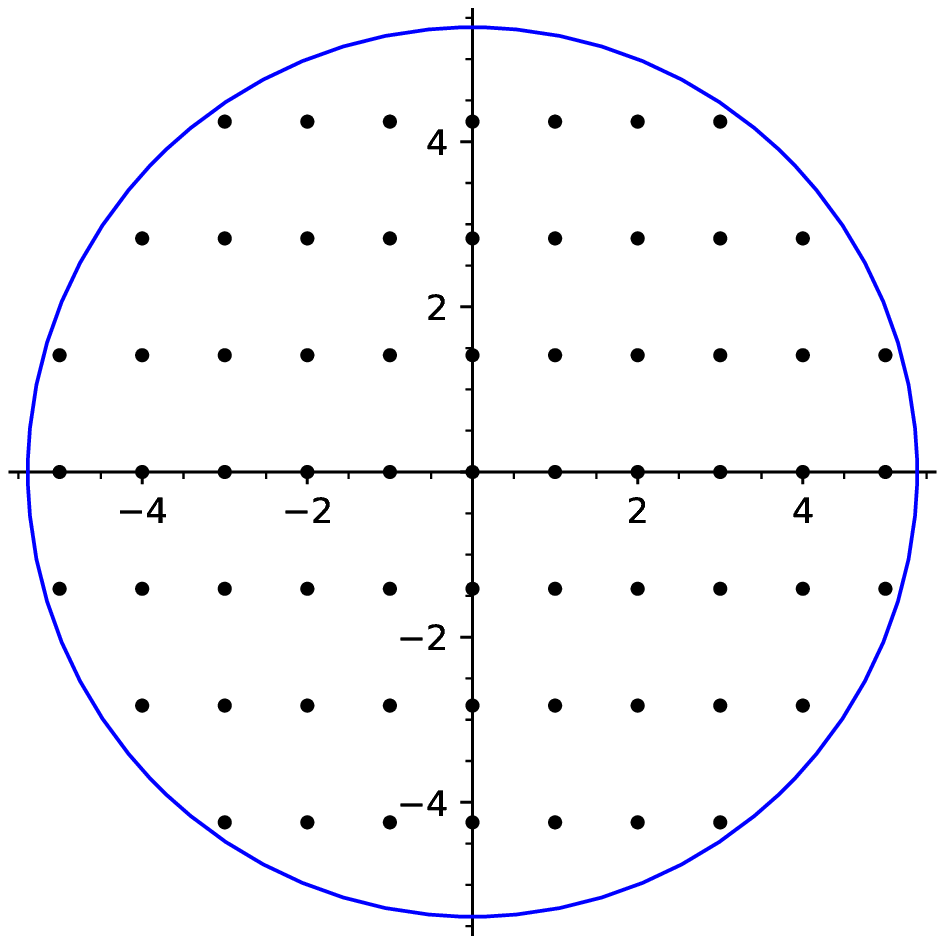}
    \caption{The 65 lattice points within the ellipse $x^2+2y^2=29$, and the 65 elements of $\Z[\sqrt{-2}]$ with norm at most 29.}
    \label{fig:x^2+2y^2=29}
\end{figure}

\begin{example}\label{ex:x^2+2y^2=29}
    We'll compute the number of lattice points contained in the ellipse defined by $x^2+2y^2=29$. We have $n=29$, $\Ceil{29/4}=8$, and $\{-29/4\}=3/4$.
    
    To compute $\Num_{8,3/4,3/4,1,2}$, we consider the sequence 
    \[
        \Floor{\frac{8-3/4}{1}+3/4},\Floor{\frac{8-3/4}{2}+3/4},\dots,\Floor{\frac{8-3/4}{8}+3/4} 
        = 8,4,3,2,2,1,1,1.
    \]
    There are 4 odd numbers, so $\Num_{8,3/4,3/4,1,2}=4$.
    
    To compute $\Num_{8,1/4,3/4,1,2}$, we consider the sequence 
    \[
        \Floor{\frac{8-3/4}{1}+1/4},\Floor{\frac{8-3/4}{2}+1/4},\dots,\Floor{\frac{8-3/4}{8}+1/4} 
        = 7,3,2,2,1,1,1,1.
    \]
    There are 6 odd numbers, so $\Num_{8,1/4,3/4,1,2}=6$.
    
    Thus, the number of lattice points in this ellipse is 
    \[
        F(29) 
        = 1+2\Num_{8,3/4,3/4,1,2}+2\Num_{8,1/4,3/4,1,2}+29+2\Floor{\frac{29}{3}}-3 
        = 1+2\cdot4+2\cdot6+29+2\cdot9-3 
        = 65.
    \]
    See Figure~\ref{fig:x^2+2y^2=29}.
\end{example}

\begin{remark}
    Consider the ring $\Z[\sqrt{-2}]$. For $z=a+b\sqrt{-2}\in\Z[\sqrt{-2}]$, the norm of $z$ is $N(z)=a^2+2b^2$. Thus, Corollary~\ref{cor:x^2+2y^2} gives a formula for the number of elements of $\Z[\sqrt{-2}]$ with norm at most $n$. Following from Example~\ref{ex:x^2+2y^2=29}, we know there are 65 elements of $\Z[\sqrt{-2}]$ with norm at most 29. We visualize these elements in Figure~\ref{fig:x^2+2y^2=29}.
\end{remark}

\bibliography{refs}

\begin{thebibliography}{1}

\bibitem{Apostol1976}
Tom~M. Apostol.
\newblock {\em Introduction to analytic number theory}.
\newblock Springer-Verlag, New York-Heidelberg, 1976.
\newblock Undergraduate Texts in Mathematics.

\bibitem{Dickson1958}
L.~E. Dickson.
\newblock {\em Introduction to the Theory of Numbers}.
\newblock Dover, 1957.

\bibitem{Grosswald1985}
Emil Grosswald.
\newblock {\em Representations of integers as sums of squares}.
\newblock Springer-Verlag, New York, 1985.

\bibitem{HardyWright1979}
G.~H. Hardy and E.~M. Wright.
\newblock {\em An introduction to the theory of numbers}.
\newblock The Clarendon Press, Oxford University Press, New York, fifth
  edition, 1979.

\bibitem{MSE-floor-summation}
Eric~Naslund. 
\newblock Rounding is asymptotically useless?
\newblock Mathematics Stack Exchange.
\newblock https://math.stackexchange.com/q/116687 (version: 2017-04-13).

\bibitem{Rademacher77}
Hans Rademacher.
\newblock {\em Lectures on elementary number theory}.
\newblock Robert E. Krieger Publishing Co., Huntington, N.Y., 1977.
\newblock Reprint of the 1964 original.

\bibitem{3b1b-pi-prime-regularities}
Grant Sanderson.
\newblock {P}i hiding in prime regularities.
\newblock 3{B}lue1{B}rown, YouTube, https://youtu.be/NaL\_Cb42WyY.

\bibitem{OEIS}
Neil J.~A. Sloane and The OEIS~Foundation Inc.
\newblock The {O}n-{L}ine {E}ncyclopedia of {I}nteger {S}equences, 2022.
\newblock http://oeis.org/.

\bibitem{sage}
W.\thinspace{}A. Stein et~al.
\newblock {\em {S}age {M}athematics {S}oftware ({V}ersion 9.5)}.
\newblock The Sage Development Team, 2022.
\newblock http://www.sagemath.org/.

\end{thebibliography}
\bibliographystyle{plain}
\end{document}